\documentclass[12pt, leqno]{article}
\usepackage{amsfonts}
\usepackage{amsmath}
\usepackage{amssymb}
\usepackage{amsthm}
\usepackage{amscd}
\usepackage{mathrsfs}
\usepackage{color}
\usepackage{latexsym}
\usepackage[pdftex]{graphicx}
\usepackage{tikz}
\usetikzlibrary{arrows.meta}
\input cyracc.def
\newfont{\tencyr}{wncyr10}
\usepackage{pgf,pgfplots}

\usetikzlibrary{calc} 
\usetikzlibrary{intersections,calc,arrows.meta}
\usetikzlibrary{hobby}
\usetikzlibrary{decorations.markings}
\usetikzlibrary{knots}
\usetikzlibrary{3d}
\tikzset{otmm/.style={x={(-135:0.5cm)},y={(0:1cm)},z={(90:1cm)}}}
\usetikzlibrary{math}

\begin{document}

\begin{center}
 {\large \bf  On a relation between Deninger's foliated dynamical systems and  Connes-Consani's adelic spaces}
 \end{center}
 
 \vspace{0.08cm}

\begin{center}

Masanori MORISHITA

\end{center}

\begin{center}

{\it Dedicated to Yuko Morishita}

\end{center}

\vspace{0.08cm}

\footnote[0]{2010 Mathematics Subject Classification. 11M55, 11R37, 13F35, 14A15, 14F20, 14F35, 37C27, 57K10, 57M10\\
Key words: Noncommutative adelic spaces, Rational Witt space, Foliated dynamical system,  Arithmetic topology, Class field theory}

{\small {\bf Abstract}:  We give a relation between  Deninger's foliated dynamical systems associated to abelian number fields and  Connes-Consani's adelic spaces. It fits with the analogy 
between knots and primes in arithmetic topology and lights up  a geometric view of class field theory.}

\vspace{0.8cm}

\begin{center}
{\bf Introduction}
\end{center}

In 1990s, Deninger and Connes proposed new approaches to study number-theoretic zeta functions. 

Motivated by his cohomological study of motivic $L$-functions, Deninger made a conjectural program for arithmetic cohomology theory attached to arithmetic schemes and, in particular, he interpreted the explicit formula in analytic number theory as Lefschetz  trace formula in his cohomological framework ([D1], [D2], [D3]). He then noticed that the would-be phase space attached to arithmetic schemes should have intimate analogies with foliated dynamical systems ([D4], [D5], [D6]).  For the case of number rings, this picture refines the analogies between knots and primes in arithmetic topology ([M]). 
Recently, motivated by the work of Kucharczyk and Scholze [KS], Deninger constructed foliated dynamical systems attached to arithmetic schemes by introducing  the notion of rational Witt space ([D7]).   

On the other hand, motivated by his study of the Riemann hypothesis, Connes introduced the ad\`{e}le class space  by discovering that it is the  noncommutative geometric object so that the space of functions on it is dual to the Bost-Connes system ([BC]), which provides the functional analytic tool in the study of the zeros of the Riemann zeta function. In particular,  he interpreted the explicit formula in analytic number theory and the Riemann hypothesis as the Lefschetz and Selberg type trace formulas, respectively,  in his noncommutative framework ([C2], [C3]).  We note that the ad\`{e}le class space may be regarded as a space of leaves on which the id\`{e}le class group acts as a flow and we should remind that foliation theory is an origin of noncommutative geometry ([C1]). 
Since then, he and Consani  have developed new foundations in number theory by means of noncommutative geometry and topos theory (for example, [CC1], [CC2]). 
Recently, they made a beautiful discovery that the analogy between knots and primes admits a geometric realization within their adelic  spaces ([CC3]).

Although both  Deninger's foliated dynamical systems and Connes-Consani' adelic spaces have the structures of foliation and dynamical system, 
their approaches seem deeply different.  Deninger's approach is concerned with arithmetic schemes and their motivic $L$-functions, while Connes' noncommutative approach to the Riemann zeta function can be generalized to automorphic $L$-functions ([S]) but not to motivic $L$-functions. Their relation has been unknown for a long time.  \\

In this paper, we give a relation between Deninger's foliated dynamical systems attached to number rings ([D7]) and Connes-Consani's adelic spaces.  Actually Deninger constructed a map from his spaces to the ad\`{e}le ring and asked if there is a relation with Connes' theory (cf. the section 7 of [D7]).  We answer to this question by showing that  Deninger's map gives rise to  a continuous map from  his foliated dynamical systems associated to abelian extensions of $\mathbb{Q}$ to Connes-Consani's adelic spaces such that it is Galois-equivariant and flow-anti-equivariant, in particular, closed orbits attached to primes in both spaces are corresponding (cf. Theorem 3.6 below).  A key observation  we found is that the arithmetic linking homomorphism of a prime $p$ with other primes
$$ {\rm lk}_p :  p^{\hat{\mathbb{Z}}} \longrightarrow \hat{\mathbb{Z}}_{(p)}^{\times} := \prod_{q \neq p} \mathbb{Z}_q ^{\times},$$
plays roles in both Deninger's theory and Connes-Consani's theory as the monodromies in the coverings of phase spaces associated to abelian extensions of $\mathbb{Q}$. In this sense, arithmetic topology provides a bridge between them, which highlights  the geometric view of class field theory.  \\

This paper is organized as follows. In Section 1, we recall Connes-Consani's adelic spaces associated to abelian extensions of $\mathbb{Q}$ and  their interpretation of class field theory within their noncommutative adelic spaces. In Section 2, we recall Deninger's foliated dynamical systems for the case of the rings of integers of number fields. Deninger's theory is about any arithmetic scheme, but, we focus on  his phase spaces attached to number rings.  In the subsection 2.3, we describe the monodromies in the coverings of Deninger systems associated to abelian extensions of $\mathbb{Q}$. In Section 3, we show a relation between Deninger's foliated dynamical system and Connes-Consani's adelic spaces associated to abelian extensions of $\mathbb{Q}$.\\
\\
{\it Acknowledgement.} I would like to thank Christopher Deninger for discussions at M\"{u}nster on August of 2023, January of 2024 and at Fukuoka on April of 2024. I would like to thank James Borger for discussions at Canberra on March of 2024. I would like to thank Alain Connes and Caterina Consani for discussions at Bologna on June of 2024 and helpful e-mail communications.  I would like to thank Jes\'{u}s Antonio \'{A}lvarez L\'{o}pez for his comment. Finally I would like to thank Jack Morava  for his encouragement and advice. The author is partially supported by the Japan Society for the Promotion of Science, KAKENHI Grant Number (C) JP22K03270.  \\
\\
{\bf Notations.}  Let $\overline{\mathbb{Q}}$ denote the algebraic closure of the rationals $\mathbb{Q}$ in the field $\mathbb{C}$ of complex numbers and let $\overline{\mathbb{Z}}$ denote the integral closure of the ring $\mathbb{Z}$ of rational integers in $\overline{\mathbb{Q}}$.  A number field means a subextension of $\overline{\mathbb{Q}}$ over $\mathbb{Q}$. For a number field $K$, ${\cal O}_K$ denotes the ring of integers in $K$, namely, ${\cal O}_K := K \cap \overline{\mathbb{Z}}$. 
Let $\mathbb{Q}^{\rm ab}$ be the maximal abelian extension of $\mathbb{Q}$ in $\overline{\mathbb{Q}}$ and let $\mathbb{Z}^{\rm ab} := {\cal O}_{\mathbb{Q}^{\rm ab}}$. 

Let $\mathbb{N}$ denote the multiplicative monoid of positive integers.  For $n \in \mathbb{N}$, let $\mu_n$ be the group of $n$-th roots of unity in $\overline{\mathbb{Z}}$, and set $\mu_{\infty} := \bigcup_{n \in \mathbb{N}} \mu_n$. For a subring $R$ of  $\overline{\mathbb{Q}}$ and $n \in \mathbb{N} \cup \{ \infty\}$, we set  $\mu_n(R) := R \cap \mu_n$.   Let $\mathbb{Q}_+$ (resp. $\mathbb{R}_+$) be the group of positive rational (resp. positive real) numbers. Note that $\mathbb{Z}^{\rm ab} = \mathbb{Z}[\mu_{\infty}]$. For a prime number $p$, let $\mu_{(p)}$ denote the group of prime-to-$p$ order roots of unity in $\overline{\mathbb{Z}}$ and  $\mu_{(p)}(R) := R \cap \mu_{(p)}$ for a subring $R$ of  $\overline{\mathbb{Q}}$ .

Let $\Sigma_{\mathbb{Q}}$ (resp. $\Sigma_{\mathbb{Q}}^{\rm f}$) denote the set of all places (resp. finite primes)  of $\mathbb{Q}$. Let $\hat{\mathbb{Z}}$ be the profinite completion of $\mathbb{Z}$, which is identified with $\prod_{p  \in \Sigma_{\mathbb{Q}}^{\rm f}}\mathbb{Z}_p$, where $\mathbb{Z}_p$ denotes the ring of  $p$-adic integers. For a prime number $p$, $\mathbb{Z}_{(p)}$ denotes the localization of $\mathbb{Z}$ at the prime ideal $(p)$ and $\hat{\mathbb{Z}}_{(p)} := \prod_{q \in \Sigma_{\mathbb{Q}}^{\rm f} \setminus \{ p \}} \mathbb{Z}_q$. The infinite prime of $\mathbb{Q}$ is denoted by $\infty$. Let $\mathbb{A}$ (resp. $\mathbb{A}^{\rm f}$) denotes the ring of ad\`{e}les (resp. finite ad\`{e}les) of $\mathbb{Q}$ so that $\mathbb{A} = \mathbb{A}^{\rm f} \times \mathbb{R}$ for the field $\mathbb{R}$ of real numbers. 
 
Let  {\bf Sch},  {\bf CRing} and {\bf Gr} denote  the categories of schemes, unital commutative  rings and groups, respectively. For a scheme $X$ with the structure sheaf ${\cal O}_X$ and a point $x \in X$,  $\kappa(x)$ denotes the residue field of the local ring ${\cal O}_{X,x}$.

For a group $Z$, a right $Z$-set $A$ and a left $Z$-set $B$, we use $A \times_{Z} B$ in two senses. Firstly $A \times_Z B$ denotes the quotient of $A \times B$ by the the left twisted $Z$-action defined by $z.(a,b) := (az^{-1}, zb)$, and secondly it denotes the quotient of  $A \times B$ by the the right twisted $Z$-action defined by $(a,b).z := (az, z^{-1} b)$, where $a \in A, b \in B$ and $z \in Z$. In this paper, we are concerned with the left twisted quotients for the cases where  
$A$ is a multiplicative abelian group $A$, $B = \mathbb{Q}_+ \, \mbox{or}\, \mathbb{R}_+$ and $Z = z^{\mathbb{Z}}$ is  a cyclic group   generated by $z$ such that there are homomorphisms $u : Z \rightarrow A$ and $v : Z \rightarrow B$. The action of $Z$ on $A$ and $B$ are multiplications via $u$ and $v$, respectively. Then $A \times_Z B$ is identified with the quotient
$Z\backslash (A \times B)$ of $A \times B$ by the diagonal action of $Z$ through the correspondence $Z\backslash(A \times B) \ni [(a,b)] \mapsto [(a^{-1},b)] \in A \times_Z B$.  We call $A \times_{z^{\mathbb{Z}}} \mathbb{R}_+$ the {\it mapping torus} of the multiplication by $u(z)$ in $A$. We are also concerned with the right twisted quotient for the case where $Z = \mathbb{Q}_+$, $A$ is a right $Z$-set, and $B = \mathbb{Q}_+ \, \mbox{or}\, \mathbb{R}_+$. Then $A \times_Z B$ is identified with the quotient $ (A \times B)/Z$ of $A \times B$ by the diagonal action of $Z$ through the correspondence $(A \times B)/Z \ni [(a,b)] \mapsto [(a,b^{-1})] \in A \times_Z B$. \\

\begin{center}
{\bf 1.  Connes-Consani's adelic spaces}
\end{center}

In this section, we recall Connes-Consani's adelic spaces associated to abelian extensions of $\mathbb{Q}$ and their geometric interpretation of class field theory  within their noncommutative adelic  spaces. For the details, readers should consult [CC3].
\\
\\
{\bf 1.1 Noncommutative adelic spaces.}  As in Notations, let $\mathbb{Q}^{\rm ab}$ be the maximal abelian extension of $\mathbb{Q}$ in $\overline{\mathbb{Q}}$, which coincides with $\mathbb{Q}(\mu_{\infty})$ by Kronecker-Weber, and let $\mathbb{Z}^{\rm ab} := {\cal O}_{\mathbb{Q}^{\rm ab}} = \mathbb{Z}[\mu_{\infty}]$.  
We identify the Galois group ${\rm Gal}(\mathbb{Q}^{\rm ab}/\mathbb{Q})$ with $\hat{\mathbb{Z}}^{\times}$ by $a \in \hat{\mathbb{Z}} \mapsto \sigma_a \in {\rm Gal}(\mathbb{Q}^{\rm ab}/\mathbb{Q})$, where $\sigma_a(\zeta) := \zeta^a$ for $\zeta  \in \mu_{\infty}$. 
Let $\mathbb{A}$ be the ring of ad\`{e}les of $\mathbb{Q}$, on which $\mathbb{Q}^{\times}$ acts by diagonal multiplication  and $\hat{\mathbb{Z}}^{\times} = \prod_p \mathbb{Z}_p^{\times}$ acts by multiplication on component-wise over all finite primes (prime numbers).  
We set
$$\mathscr{X}_{\mathbb{Q}} := \mathbb{Q}^{\times}\backslash \mathbb{A}/\hat{\mathbb{Z}}^{\times}, \; \mbox{and}\; \mathscr{X}_{\mathbb{Q}^{\rm ab}} :=
\mathbb{Q}^{\times} \backslash \mathbb{A},     \leqno{(1.1.1)}$$ 
where $\mathbb{Q}^{\times}$ (resp. $\hat{\mathbb{Z}}^{\times}$) acts on $\mathbb{A}$ through the diagobal embedding $\mathbb{Q} \hookrightarrow \mathbb{A}$ (resp. the inclusion $\hat{\mathbb{Z}} \subset \mathbb{A}^{\rm f}$).  We regard $\mathscr{X}_{\mathbb{Q}}, \mathscr{X}_{\mathbb{Q}^{\rm ab}}$ as the noncommutative quotients of $\mathbb{A}$ corresponding to crossed product noncommutative algebras. The noncommutative space $\mathscr{X}_{\mathbb{Q}^{\rm ab}}$ was introduced by Connes to study the zeros of the Riemann zeta function ([C1]) and so we call $\mathscr{X}_{\mathbb{Q}^{\rm ab}}$ the {\it Connes ad\`{e}le class space}. In the following,  $\mathscr{X}_{\mathbb{Q}}$ and $\mathscr{X}_{\mathbb{Q}^{\rm ab}}$ are seen as noncommutative geometric models of ${\rm Spec}(\mathbb{Z})$ and ${\rm Spec}(\mathbb{Z}^{\rm ab})$, respectively. We also note that $\mathscr{X}_{\mathbb{Q}}$ corresponds to the space of points of the associated topos, called the {\it scaling topos} (cf. [CC1], [CC2]).
  
 More generally, we also consider the noncommutative spaces for all abelian extensions of $\mathbb{Q}$, namely, subextensions of $\mathbb{Q}^{\rm ab}$.  Note that an abelian extension $K$ of $\mathbb{Q}$ corresponds bijectively to a closed subgroup  $U_K := {\rm Gal}(\mathbb{Q}^{\rm ab}/K)$ of ${\rm Gal}(\mathbb{Q}^{\rm ab}/\mathbb{Q}) = \hat{\mathbb{Z}}^{\times}$.  Note that $U_{\mathbb{Q}} = \hat{\mathbb{Z}}^{\times}$ and $U_{\mathbb{Q}^{\rm ab}} = \{ 1 \}$. We then define the noncommutative space $\mathscr{X}_K$  by
$$ \mathscr{X}_K := \mathbb{Q}^{\times} \backslash \mathbb{A} /U_K,  \leqno{(1.1.2)}$$
which we call the {\it Connes-Consani adelic space} associated to an abelian extension $K$ of $\mathbb{Q}$.  We note that the Galois group
${\rm Gal}(K/\mathbb{Q})$ acts on $\mathscr{X}_K$ from the right by
$$ \mathscr{X}_K \times {\rm Gal}(K/\mathbb{Q}) \longrightarrow \mathscr{X}_K; \;\; (\mathbb{Q}^{\times}\alpha U_K, g) \mapsto \mathbb{Q}^{\times}\alpha \tilde{g}U_K, \leqno{(1.1.3)}$$
where $\tilde{g}$ is an extension of $g$ to ${\rm Gal}(\mathbb{Q}^{\rm ab}/\mathbb{Q}) = \hat{\mathbb{Z}}^{\times}$. For another extension
$\tilde{g}'$ of $g$, we may write $\tilde{g}' = \tilde{g} u$ for some $u \in U_K$ and so the action is well defined.

The Connes ad\`{e}le class space $\mathscr{X}_{\mathbb{Q}^{\rm ab}}$ in (1.1.1) is equipped with the dynamical system defined by the multiplicative action of the id\`{e}le class  group, and the Connes-Consani adelic space $\mathscr{X}_K$ is equipped with dynamical system defined by the action of $\mathbb{R}_+$ on the $\infty$-component of $\mathbb{A}$ by multiplication:
$$ \mathscr{X}_K \times \mathbb{R}_+ \longrightarrow \mathscr{X}_K; \;\; (\mathbb{Q}^{\times}(a_v)U_K, t) \mapsto   \mathbb{Q}^{\times}((a_q), a_{\infty}t)U_K. \leqno{(1.1.4)}$$
Since $\mathscr{X}_K$ is the quotient space, it has the foliation structure taking the cosets as leaves. The Connes-Consani adelic space  $\mathscr{X}_K$ may be regarded as a  noncommutative geometric model for ${\rm Spec}({\cal O}_K)$. We let   $\varphi_{\mathbb{Q}^{\rm ab}/K}  : \mathscr{X}_{\mathbb{Q}}^{\rm ab} \rightarrow \mathscr{X}_{K}$ and $\varphi_{K/\mathbb{Q}} : \mathscr{X}_{K} \rightarrow \mathscr{X}_{\mathbb{Q}}$ be the natural coverings.

In the Connes-Consani adelic space $\mathscr{X}_K$ with $K$ being of finite degree over $\mathbb{Q}$, a non-zero prime $\frak{p} \in {\rm Spec}({\cal O}_K)$ whose residue field is of characteristic $p$ is visualized as the $\mathbb{R}_+$-orbit $C_{\frak{p}}$ of the ad\`{e}le $(a_v)$ with $a_p=0$, $a_v = 1$ ($v \neq p, \infty$) and $a_{\infty} > 0$:
$$ C_{\frak{p}} := \mathbb{Q}^{\times} \{ (a_v) \in \mathbb{A} \; \mid \; a_p = 0, a_q = 1 (q \neq p), a_{\infty} > 0 \} U_K.
 \leqno{(1.1.5)}$$
By the projection to the $\infty$-component, $C_{\frak{p}}$ is isomorphic to the circle $\mathbb{R}_+/{\rm N}\frak{p}^{\mathbb{Z}} \simeq \mathbb{R}/(\log {\rm N}\frak{p})\mathbb{Z}$ of length $\log {\rm N}\frak{p}$.\\
\\
{\bf 1.2   Abelian coverings of adelic spaces and class field theory.}  From the viewpoint of arithmetic topology, class field theory describes how primes are linked each other. For instance, the mod 2 {\it arithmetic linking number} of an odd prime $p$ with another odd prime $q$ is nothing but the Legendre symbol $\displaystyle{ \left(   \frac{q}{p} \right) }$, which is described as the image of $p$ in $\mathbb{F}_q^{\times}/(\mathbb{F}_q^{\times})^2 = \{ \pm 1 \}$ (See [M; Chapter 4]). So the arithmetic linking of a prime $p$ with all other primes is described by the image $p$ under the following diagonal homomorphism, which we call  the {\it arithmetic linking homomorphism for a  prime $p$}:
$$ {\rm lk}_p \; : \; p^{\hat{\mathbb{Z}}} \longrightarrow \hat{\mathbb{Z}}_{(p)}^{\times} = \prod_{q \neq p} \mathbb{Z}_{q}^{\times}. \leqno{(1.2.1)}$$
 Connes and Consani noticed that the arithmetic linking homomorphism coincides with the following homomorphism of \'{e}tale fundamental groups, which we denote by the same ${\rm lk}_p$, 
$${\rm lk}_p \; : \; \pi_1({\rm Spec}(\mathbb{F}_p)) \longrightarrow \pi_1^{\rm ab}({\rm Spec}(\mathbb{Z}_{(p)})) \leqno{(1.2.2)}$$
induced by the reduction map  $\mathbb{Z}_{(p)} \rightarrow \mathbb{F}_p$ modulo $p$, where $\pi_1^{\rm ab}$ denotes the abelianized \'{e}tale fundamental group.  Here the prime $p$ corresponds to the Frobenius automorphism $\sigma_p \in \pi_1({\rm Spec}(\mathbb{F}_p))$
 and so we identify ${\rm lk}_p(p)$ with ${\rm lk}_p(\sigma_p)$. 
In terms of field extension,  $\pi_1({\rm Spec}(\mathbb{F}_p)) = {\rm Gal}(\overline{\mathbb{F}}_p/\mathbb{F}_p)$, where $\overline{\mathbb{F}}_p$ is an algebraic closure of $\mathbb{F}_p$, and $\pi_1^{\rm ab}({\rm Spec}(\mathbb{Z}_{(p)}))$ is the Galois group of the maximal abelian extension $\mathbb{Q}_{(p)}^{\rm ab}$ of $\mathbb{Q}$ such that $p$ is unramified. By Kronecker-Weber, 
$\mathbb{Q}_{(p)}^{\rm ab}$ coincides with $\mathbb{Q}(\mu_{(p)})$ (cf. Notations for $\mu_{(p)}$) and we have the canonical identification
$$ \hat{\mathbb{Z}}_{(p)}^{\times} \stackrel{\sim}{\longrightarrow} {\rm Gal}(\mathbb{Q}(\mu_{(p)})/\mathbb{Q}); \;\; a \mapsto (\mu_{(p)} \ni \zeta \mapsto \zeta^a \in \mu_{(p)}).\;\; $$
Then the linking homomorphism in (1.2.1) and (1.2.2)  is also identified with 
the homomorphism
$$ {\rm lk}_p :  {\rm Gal}(\overline{\mathbb{F}}_p/\mathbb{F}_p) \longrightarrow {\rm Gal}(\mathbb{Q}(\mu_{(p)})/\mathbb{Q}) \leqno{(1.2.3)} $$
induced by the reduction map $\mathbb{Z}[\mu_{(p)}] \rightarrow \overline{\mathbb{F}}_p$ modulo a maximal ideal of $\mathbb{Z}[\mu_{(p)}]$ over
$p$, where ${\rm lk}_p(\sigma_p)$ is the Frobenius automorphism over $p$ in the extension $\mathbb{Q}(\mu_{(p)})/\mathbb{Q}$ given by
$\zeta \mapsto \zeta^p$ for $\zeta \in \mu_{(p)}$.

Let $\eta : \hat{\mathbb{Z}}_{(p)}^{\times} \times \mathbb{R}_+ \rightarrow \mathbb{A}$ be the map defined by 
$$\eta((a_q)_{q \neq p}, s)) := (a_v), \;\; a_v :=  \left\{ \begin{array}{ll} 0 \; & \; {\it if}\; v = p,  \\
 a_q  \; &\;  {\it if} \; v = q \neq p, \\
 s  \;  & \; {\it if}\; v = \infty.  \end{array} \right.$$
Then Connes and Consani have made the following beautiful discovery that the arithmetic linking homomorphism ${\rm lk}_p$ gives really a geometric  monodromy around the $\mathbb{R}_+$-orbit $C_p$ in (1.1.5) in the natural covering $\varphi_{\mathbb{Q}^{\rm ab}/\mathbb{Q}} : \mathscr{X}_{\mathbb{Q}^{\rm ab}} \rightarrow \mathscr{X}_{\mathbb{Q}}$. We give their result ([CC3; Theorem 0.2, Theorem 3.2]) in a slightly generalized form. The proof is similar to those of [ibid]. \\
\\
{\bf Theorem 1.2.4} (cf. [CC3; Theorem 0.2, Theorem 3.2]). {\it Notations being as above, let $L$ be a number field containing $\mathbb{Q}(\mu_{(p)})$.}  {\it Then the  map $\eta$ induces a canonical isomorphism of $\varphi_{L/\mathbb{Q}}^{-1}(C_p)$ with the mapping torus of the multiplication by the ${\rm lk}_p(p) (= {\rm lk}_p(\sigma_p))$ in $\hat{\mathbb{Z}}_{(p)}^{\times} (= \pi_1^{\rm ab}({\rm Spec}(\mathbb{Z}_{(p)}))$}\footnote{In [CC3],  the quotient $p^{\mathbb{Z}}\backslash (\hat{\mathbb{Z}}_{(p)}^{\times} \times \mathbb{R}_+)$ by the diagonal action of $p^{\mathbb{Z}}$ is considered. As is noted in Notations,  $p^{\mathbb{Z}}\backslash (\hat{\mathbb{Z}}_{(p)}^{\times} \times \mathbb{R}_+)$ is identified with $\hat{\mathbb{Z}}_{(p)}^{\times} \times_{p^{\mathbb{Z}}} \mathbb{R}_+$ through the correspondence $[(a,b)] \mapsto [(a^{-1},b)]$.}:
$$ \varphi_{L/\mathbb{Q}}^{-1} (C_p) \simeq  \hat{\mathbb{Z}}_{(p)}^{\times} \times_{p^{\mathbb{Z}}} \mathbb{R}_+.$$
{\it The monodromy  around the circle $C_p$ in  $\varphi_{L/\mathbb{Q}}^{-1}(C_p)$ is given by the multiplication by ${\rm lk}_p(p) (= {\rm lk}_p(\sigma_p))$ in $\hat{\mathbb{Z}}_{(p)}^{\times} (= \pi_1^{\rm ab}({\rm Spec}(\mathbb{Z}_{(p)}))$.}\\

Next we shall present the results similar to Theorem 1.2.4  for coverings of  Connes-Consani adelic spaces attached to finite abelain extensions of $\mathbb{Q}$. Let $F$ be a finite abelian subextension of $\mathbb{Q}^{\rm ab}/\mathbb{Q}$ and let 
$$\chi_F : \hat{\mathbb{Z}}^{\times} = {\rm Gal}(\mathbb{Q}^{\rm ab}/\mathbb{Q}) \longrightarrow {\rm Gal}(F/\mathbb{Q})$$
 be the restriction homomorphism. Let $R_F$ be the set of finite primes of $\mathbb{Q}$ which are ramified in $F$, namely, 
$$R_F := \{ p \in \Sigma_{\mathbb{Q}}^{\rm f} \; | \; \chi(\mathbb{Z}_p^{\times}) \neq \{ 1 \} \}.$$
Then $\chi_F$ factors through $\prod_{q \in R_F} \mathbb{Z}_q^{\times}$:
$$ \chi_F : \prod_{q \in R_F} \mathbb{Z}_q^{\times}  \rightarrow {\rm Gal}(F/\mathbb{Q}). $$
 Let $p \in \Sigma_{\mathbb{Q}}^{\rm f} \setminus R_F$, namely,  suppose that $p$ is unramified in $F/\mathbb{Q}$, equivalently, $F \subset \mathbb{Q}(\mu_{(p)})$. 
 Viewing $p \in \prod_{q \in R_F} \mathbb{Z}_q^{\times}$, $\chi_F (p)$ coincides with the Artin symbol over $p$ in $F/\mathbb{Q}$:
$$ \chi_F(p) = \left(  \frac{F/\mathbb{Q}}{p} \right) \in {\rm Gal}(F/\mathbb{Q}). \leqno{(1.2.5)}$$
Let $\mathscr{X}_F$ be the Connes-Consani adelic space  for $K = F$ in  (1.1.2) and let $\varphi_{F/\mathbb{Q}} : \mathscr{X}_F \rightarrow \mathscr{X}_{\mathbb{Q}}$ be the natural covering.  For $x = [(a_v)] \in \mathscr{X}_{\mathbb{Q}}$, we set
$$ Z(x) := \{ v \in \Sigma_{\mathbb{Q}} \; | \; a_v = 0 \}.$$
For a finite subset $S$ of $\Sigma_{\mathbb{Q}}$, we say that $\varphi_{F/\mathbb{Q}}$ is {\it unramified outside} $S$ if the following condition is satisfied for any $x \in \mathscr{X}_{\mathbb{Q}}$:
$$ Z(x) \cap S = \emptyset \Longrightarrow {\rm Gal}(F/\mathbb{Q}) \; \mbox{acts freely transitively on}\; \varphi_{F/\mathbb{Q}}^{-1}(x).$$
Then we see the following\\
\\
{\bf Proposition 1.2.6} ([CC3; Proposition 2.6]). {\it Notation being as above, $\varphi_{F/\mathbb{Q}} : \mathscr{X}_F \rightarrow \mathscr{X}_{\mathbb{Q}}$ is unramified outside $R_F \cup \{ \infty \}$.}\\
\\
Now we have the following theorem, which is similar to Theorem 1.2.3 but more precise, for the covering $\varphi_{F/\mathbb{Q}} : \mathscr{X}_F \rightarrow \mathscr{X}_{\mathbb{Q}}$.\\
\\
{\bf Theorem 1.2.7} ([CC3; Theorem 0.1, Theorem 2.9]).  {\it Notations being as above, let $p \in \Sigma_{\mathbb{Q}}^{\rm f} \setminus R_F$.
Then the following assertions hold.}\\
(1) {\it The map $\eta$ in Theorem 1.2.3 induces a canonical isomorphism of $\varphi_{F/\mathbb{Q}}^{-1}(C_p)$ with the mapping torus of the multiplication by $\displaystyle{\chi_F(p) = \left(  \frac{F/\mathbb{Q}}{p} \right)}$ in ${\rm Gal}(F/\mathbb{Q})$}:
$$\varphi_{F/\mathbb{Q}}^{-1} (C_p)  \simeq  {\rm Gal}(F/\mathbb{Q}) \times_{p^{\mathbb{Z}}} \mathbb{R}_+.$$
{\it The monodromy around the circle  $C_p$ in $\varphi_{F/\mathbb{Q}}^{-1}(C_p)$  is given by the multiplication by the Artin symbol $\displaystyle{\chi_F(p) = \left(  \frac{F/\mathbb{Q}}{p} \right)}$ in ${\rm Gal}(F/\mathbb{Q})$. }\\
(2) {\it Let $(p) = \frak{p}_1 \cdots \frak{p}_r$ be the decomposition of $(p)$ into distinct primes in $F$ such that $f = {\rm deg}(\frak{p}_i) = [\kappa(\frak{p}_i) : \mathbb{F}_p]$ and $fr = [F:\mathbb{Q}]$. Then $\varphi_{F/\mathbb{Q}}^{-1}(C_p)$ is decomposed into connected components $($circles$)$ in the same manner, namely, 
$\varphi_{F/\mathbb{Q}}^{-1}(C_p) = \tilde{C}_1 \sqcup \cdots \sqcup \tilde{C}_r$
 such that $\varphi_{F/\mathbb{Q}}|_{\tilde{C}_i} : \tilde{C}_i \rightarrow C_p$ is a cyclic covering of degree $f$ and $fr = [F:\mathbb{Q}]$.} \\
\\
{\bf Example 1.2.8.} Let $p$ and $q$ be distinct odd prime numbers and consider the case that $F = \mathbb{Q}(\sqrt{q})$. Then the monodromy around the circle
$C_p$ in the covering $\varphi_{\mathbb{Q}(\sqrt{q})/\mathbb{Q}}$ is given by the multiplication by the Legendre symbol $\displaystyle{\left( \frac{q}{p} \right)}$ in ${\rm Gal}(\mathbb{Q}(\sqrt{q})/\mathbb{Q}) = \{ \pm 1 \}$ such that we have
$$ \begin{array}{ccc}
\varphi_{\mathbb{Q}(\sqrt{q})/\mathbb{Q}}^{-1}(C_p) = \tilde{C}_1 \cup \tilde{C}_2  \; (\mbox{disjoint}) & \Longleftrightarrow & \displaystyle{\left( \frac{q}{p} \right)}= 1,\\
\varphi_{\mathbb{Q}(\sqrt{q})/\mathbb{Q}}^{-1}(C_p) = \tilde{C}   & \Longleftrightarrow & \displaystyle{\left( \frac{q}{p} \right)} = -1.\\
\end{array}
$$
\\

\begin{center}
{\bf 2.  Deninger's foliated dynamical systems for number rings}
\end{center}

In the subsections 2.1 and 2.2, we recall Deninger's construction of the foliated dynamical systems attached to number rings and his description of the structure of $\mathbb{R}_+$-orbits. For the details including the case of higher dimensional arithmetic schemes, readers should consult [D7].   In  the subsection 2.3, we describe the monodromies of coverings of certain dynamical systems associated to abelian extensions of $\mathbb{Q}$, which we call {\it Deninger systems}. \\
\\
{\bf 2.1 $\mathbb{C}$-valued points of rational Witt spaces.} For a commutative ring $R$, the ring $W_{\rm rat}(R)$ of {\it rational Witt vectors}  is defined by 
$$ W_{\rm rat}(R) := \left\{ \frac{P(t)}{Q(t)}\,  \middle| \, P(t), Q(t) \in R[t], \; P(0)=Q(0) = 1 \right\}.$$
It is a subring of the big Witt ring $W(R) := 1 + tR[[t]]$ equipped with the addition
$f \oplus g := f \cdot g$ for $f, g \in W(R)$ and the multiplication satisfying $(1 - at)\odot (1-bt) = (1-abt)$ for $a, b \in R$.  It also carries the $n$-th Frobenius endomorphism $F_n : W_{\rm rat}(R)  \longrightarrow W_{\rm rat}(R)$
 for all $n \in \mathbb{N}$. A ring homomorphism $R \rightarrow R'$ induces a ring homomorphism $W_{\rm rat}(R) \rightarrow W_{\rm rat}(R')$ in the obvious manner and hence $W_{\rm rat}$ gives a functor on ${\bf CRing}$. There is the multiplicative injective Teichm\"{u}ller map
$$ [\; ] : R \longrightarrow W_{\rm rat}(R); \;\; a \mapsto [a] := 1 - at,$$
which splits by the canonical surjective ring homomorphism
$$ W_{\rm rat}(R) \longrightarrow R; \;\; f \mapsto - f'(0).$$
We note that $F_n([a]) = [a^n]$ for $a \in R$ and $n \in \mathbb{N}$. 
If a  group $G$ acts on a commutative ring $R$, then $G$ acts on $W_{\rm rat}(R)$ through the action of $G$ on coefficients. \\
\\
{\bf Proposition 2.1.1} ([KS; Lemma 4.9], [D7; Theorem 1.5]).  {\it Let $G$ be a profinite group acting continuously on a commutative ring $R$. Then the natural inclusion $W_{\rm rat}(R^G) \hookrightarrow W_{\rm rat}(R)^G$ induces an isomorphism} 
$$ W_{\rm rat}(R^G) \stackrel{\sim}{\longrightarrow} W_{\rm rat}(R)^G. $$
 {\it For an algebraically closed field $\Omega$, we have the following bijection}
$$ {\rm Hom}_{\bf CRing}(W_{\rm rat}(R), \Omega)/G \stackrel{\sim}{\longrightarrow} {\rm Hom}_{\bf CRing}(W_{\rm rat}(R^G), \Omega), \;\; \alpha \mapsto \alpha|_{W_{\rm rat}(R^G)}, $$
{\it where $G$ acts on the set ${\rm Hom}_{\bf CRing}(W_{\rm rat}(R),\Omega)$ by $\alpha \mapsto \alpha \circ \sigma$ for $\sigma \in G$.}
\\ 
\\
{\bf Example 2.1.2} ([KS; Remark 4.3], [D7; Proposition 1.4]).  For an algebraically closed field $\Omega$, there is the ring isomorphism
$$  \mathbb{Z}[\Omega^{\times}] \stackrel{\sim}{\longrightarrow} W_{\rm rat}(\Omega); \;\; \sum_a {n_a} a \mapsto \prod_a (1-at)^{n_a}, $$
 where $\mathbb{Z}[\Omega^{\times}]$ denotes the group ring of the multiplicative (commutative) group $\Omega^{\times}$ over $\mathbb{Z}$. When $\Omega$ is an algebraic closure $\overline{k}$ of a perfect field $k$, then the above isomorphism commutes with the actions of the Galois group ${\rm Gal}(\overline{k}/k)$ on the both sides. By Proposition 2.1.1, the above isomorphism induces the ring isomorphism
 $$ \mathbb{Z}[\overline{k}^{\times}]^{{\rm Gal}(\overline{k}/k)} \stackrel{\sim}{\longrightarrow} W_{\rm rat}(k).$$
 \vspace{0.02cm}
 
For  a scheme $X$ with the structure sheaf ${\cal O}_X$, the {\it rational Witt space} of $X$ is defined by the ringed  space
$$ W_{\rm rat}(X) := (X_{\rm top}, W_{\rm rat}({\cal O}_X)).$$
Here $X_{\rm top}$ denotes the underlying topological space and $W_{\rm rat}({\cal O}_X)$ be the sheafification of the Zariski presheaf associating $W_{\rm rat}({\cal O}_X(U))$ to  an open subset $U$ of $X$. For a point $x \in X$, we have $W_{\rm rat}({\cal O}_X)_x = W_{\rm rat}({\cal O}_{X,x})$. For $n \in \mathbb{N}$,  the $n$-th Frobenius endomorphism $F_n : W_{\rm rat}(X) \rightarrow W_{\rm rat}(X)$ is defined by the identity on $X_{\rm top}$ and $F_n$ on $W_{\rm rat}({\cal O}_X)$ induced by $F_n : W_{\rm rat}({\cal O}_X(U)) \rightarrow W_{\rm rat}({\cal O}_X(U))$ for open subsets $U$ of $X$.   

For  schemes $S$ and $X$, a morphism $f : S \rightarrow W_{\rm rat}(X)$ is defined by a morphism of ringed spaces $(f_{\rm top}, f^{\#}) : (S_{\rm top}, {\cal O}_S) \rightarrow (X_{\rm top}, W_{\rm rat}({\cal O}_X))$, which is local in the following sense: For $s \in S$, there is a unique ring homomorphism $\tilde{f}^{\#}_s : W_{\rm rat}(\kappa(f(s))) \rightarrow \kappa(s)$
such that the following  diagram commutes.
$$  \begin{array}{ccc}
W_{\rm rat}({\cal O}_{X, f(s)}) & \stackrel{ f^{\#}_s}{\longrightarrow} & {\cal O}_S \\
\downarrow & & \downarrow \\
W_{\rm rat}(\kappa(f(s))) & \stackrel{ \tilde{f}^{\#}_s}{\longrightarrow} & \kappa(s).
\end{array}
$$
We define the set  $W_{\rm rat}(X)(S)$ of $S$-valued points of $W_{\rm rat}(X)$ by the set of morphisms $S \rightarrow W_{\rm rat}(X)$. 
For a field $k$, we set $W_{\rm rat}(X)(k) := W_{\rm rat}(X)({\rm Spec}(k))$, called the set of $k$-valued points of $W_{\rm rat}(X)$.  
By definition, $W_{\rm rat}(X)(k)$ is given by 
$$  W_{\rm rat}(X)(k) = \{ (x, P)  \mid  x \in X, \; P \in {\rm Hom}_{\bf CRing}(W_{\rm rat}(\kappa(x)), k) \}. $$
The canonical surjective ring homomorphisms $W_{\rm rat}({\cal O}_X(U)) \rightarrow {\cal O}_X(U)$ for open subsets $U$ of $X$ induce a surjective morphism of sheaves  of commutative rings $ f^{\#} : W_{\rm rat}({\cal O}_X) \rightarrow {\cal O}_X$. This gives a morphism
$$ f = ({\rm id}, f^{\#}) : X \longrightarrow W_{\rm rat}(X).$$
Then  we have a canonical injective map 
$$ X(k) := {\rm Hom}_{\bf Sch}({\rm Spec}(k), X) \stackrel{\subset}{\longrightarrow} W_{\rm rat}(X)(k); \;\; g  \mapsto f \circ g.$$
We also have a multiplicative map
$$ [\; ] : \Gamma(X, {\cal O}_X) \longrightarrow {\rm Map}(W_{\rm rat}(X)(k), \Omega);\; h \mapsto ((x, P) \mapsto P([h(x)])).$$
For example, if  $X = {\rm Spec}(\mathbb{Z})$ = the final object of {\bf Sch} and $k = \mathbb{C}$,  then ${\rm Spec}(\mathbb{Z})(\mathbb{C})$ consists of a single element, while $W_{\rm rat}({\rm Spec}(\mathbb{Z}))(\mathbb{C})$ is of infinite-dimensional over $\mathbb{C}$, and the  multiplicative map
$$ [\; ] : \mathbb{Z} \longrightarrow {\rm Map}(W_{\rm rat}({\rm Spec}(\mathbb{Z}))(\mathbb{C}), \mathbb{C}) $$
 enables us to regard  an integer $n \in \mathbb{Z}$  as a $\mathbb{C}$-valued function $[n]$ on $W_{\rm rat}({\rm Spec}(\mathbb{Z}))(\mathbb{C})$. \\

In the following, we focus on the arithmetic situation.   Let $K$ be a number field.  We set 
$$X_K := {\rm Spec}({\cal O}_K)$$
and 
$$ \dot{X}_K(\mathbb{C}) := W_{\rm rat}(X_K)(\mathbb{C}).$$
So we have 
$$\dot{X}_K(\mathbb{C}) = \{ (\frak{p}, P) \mid  \frak{p} \in X_K,  \; P \in {\rm Hom}_{\bf CRing}(W_{\rm rat}(\kappa(\frak{p})), \mathbb{C}) \}. \leqno{(2.1.3)}$$
We note that $\dot{X}_K(\mathbb{C})$ carries the Frobenius $\mathbb{N}$-action defined by
 $$ F_n(\frak{p}, P) := (\frak{p}, P \circ F_n) \;\; \mbox{for}\; n \in \; \mathbb{N}. \leqno{(2.1.4)}$$
We note by Example 2.1.2 that  when $\kappa(\frak{p})$ is algebraically closed,  giving a point $P \in {\rm Hom}_{\bf CRing}(W_{\rm rat}(\kappa(\frak{p})), \mathbb{C})$ is equivalent to giving $P^{\times} \in {\rm Hom}_{\bf Gr}(\kappa(\frak{p})^{\times}, \mathbb{C}^{\times})$. For example,  $\kappa(\frak{p})$ is algebraically closed, if $K$ contains $\mathbb{Q}(\mu_{\infty})$ and $\frak{p}$ is a closed point in $X_K$. In particular, when $K = \overline{\mathbb{Q}}$, we  have
$$ \dot{X}_{\overline{\mathbb{Q}}}(\mathbb{C}) = \{ (\frak{p}, P^{\times}) \mid  \frak{p} \in X_{\overline{\mathbb{Q}}},  \; P^{\times} \in {\rm Hom}_{\bf Gr}(\kappa(\frak{p})^{\times}, \mathbb{C}^{\times}) \}.
\leqno{(2.1.5)}
 $$
 The next proposition gives an alternative description of $\dot{X}_{\overline{\mathbb{Q}}}(\mathbb{C})$.\\
 \\
 {\bf Proposition 2.1.6} ([D7, Remark 3.4]).  {\it A point $(\frak{p}, P^{\times}) \in \dot{X}_{\overline{\mathbb{Q}}}(\mathbb{C})$ can be  identified with the multiplicative map 
 $P : \overline{\mathbb{Z}} \rightarrow \mathbb{C}$ satisfying the conditions}:  (i) $P(0) = 0, P(1)=1$, (ii) {\it $P^{-1}(0)$ is a prime ideal}, (iii) {\it $P$ factors through $\overline{\mathbb{Z}}/P^{-1}(0)$.  Indeed such a map $P$ gives $\frak{p}= P^{-1}(0)$ and induces a homomorphism $P^{\times} : \kappa(\frak{p})^{\times} \rightarrow \mathbb{C}^{\times}$. }\\
 \\
 By Proposition 2.1.6,  we can define a topology on $\dot{X}_{\overline{\mathbb{Q}}}(\mathbb{C})$. Namely, since $\dot{X}_{\overline{\mathbb{Q}}}(\mathbb{C})$ is the set of multiplicative maps, we equip it  with the topology of pointwise convergence.  \\ 
 
 Let $L$ be Galois extension of $K$ in $\overline{\mathbb{Q}}$.  The Galois group ${\rm Gal}(L/K)$ acts continuously on $\dot{X}_L(\mathbb{C})$ from the right by
 $$ (\frak{p}, P)^{\sigma} := (\frak{p}^{\sigma}, P \circ \sigma) \; \; \mbox{for}\; \sigma \in \; {\rm Gal}(L/K),  \leqno{(2.1.7)}$$
 where a point  $\frak{p}^{\sigma} \in X_L$ is nothing but $\sigma^{-1}(\frak{p})$ as a prime ideal of ${\cal O}_L$ and so
 $$ P \circ \sigma :  W_{\rm rat}(\kappa(\sigma^{-1}(\frak{p}))) \stackrel{\sigma}{\longrightarrow} W_{\rm rat}(\kappa(\frak{p})) \stackrel{P}{\longrightarrow} \mathbb{C}.$$
The action of ${\rm Gal}(L/K)$ commutes with the Frobenius $\mathbb{N}$-action. The natural morphism $\pi_{L/K} : X_L \rightarrow X_K$ 
induces the $\mathbb{N}$-equivariant map
 $$ \dot{\pi}_{L/K} : \dot{X}_L(\mathbb{C}) \longrightarrow \dot{X}_K(\mathbb{C}); (\frak{P}, P) \mapsto (\frak{p}, P|_{W_{\rm rat}(\kappa(\frak{p}))}) \;\; (\frak{p} = \frak{P} \cap K). \leqno{(2.1.8)}$$
 By Proposition 2.1.1, we obtain the following descent property.\\
\\
{\bf Proposition 2.1.9} ([D7; Corollary 3.2, Corollary 3.3]). {\it Notations being as above, $\dot{\pi}_{L/K}$ induces the following  bijection, which is equivariant under the Frobenius $\mathbb{N}$-action}
$$ \dot{X}_L(\mathbb{C})/{\rm Gal}(L/K) \stackrel{\sim}{\longrightarrow} \dot{X}_K(\mathbb{C}).$$
{\it In particular, by (2.1.5), we have the the Frobenius $\mathbb{N}$-equivariant bijections}
$$ \begin{array}{ll} \dot{X}_K(\mathbb{C}) & \simeq  \dot{X}_{\overline{\mathbb{Q}}}(\mathbb{C})/{\rm Gal}(\overline{\mathbb{Q}}/K)
 \\  & = \{ (\frak{p}, P^{\times}) \mid  \frak{p} \in X_{\overline{\mathbb{Q}}},  \; P^{\times} \in {\rm Hom}_{\bf Gr}(\kappa(\frak{p})^{\times}, \mathbb{C}^{\times}) \}/{\rm Gal}(\overline{\mathbb{Q}}/K).
\end{array}$$
\\
By Proposition 2.1.9, we equip $\dot{X}_K(\mathbb{C})  \simeq \dot{X}_{\overline{\mathbb{Q}}}(\mathbb{C})/{\rm Gal}(\overline{\mathbb{Q}}/K)$ the quotient topology so that  the map $\dot{\pi}_{L/K}$ is continuous and the bijection $ \dot{X}_L(\mathbb{C})/{\rm Gal}(L/K) \stackrel{\sim}{\rightarrow} \dot{X}_K(\mathbb{C})$ is a homeomorphism.

We order $\mathbb{N}$ by division and set
$$ \check{X}_K(\mathbb{C}) := \varinjlim_{n \in \mathbb{N}} \dot{X}_K(\mathbb{C})$$
equipped with  the inductive limit topology. 
We denote  the image of $(\frak{p}, P)$ under the inclusion $ \dot{X}_K(\mathbb{C}) \hookrightarrow \check{X}_K(\mathbb{C})$
by the same $(\frak{p}, P)$.
Then $\check{X}_K(\mathbb{C})$ carries the inverted Frobenius $\mathbb{Q}_+$-action inherited from (2.1.4):
$$ F_q((\frak{p}, P)) := (\frak{p}, P \circ F_q) \;\; \mbox{for}\; q \in \; \mathbb{Q}_+. \leqno{(2.1.10)}$$
It also carries the continuous ${\rm Gal}(L/K)$-action inherited from (2.1.7).
The map $\dot{\pi}_{L/K}$ in (2.1.8) induces the $\mathbb{Q}_+$-equivariant continuous map
$$ \check{\pi}_{L/K} : \check{X}_L(\mathbb{C}) \longrightarrow \check{X}_K(\mathbb{C}). \leqno{(2.1.11)}$$
By Proposition 2.1.9, it induces a $\mathbb{Q}_+$-equivariant homeomorphism: 
$$ \check{X}_L(\mathbb{C})/{\rm Gal}(L/K) \stackrel{\sim}{\longrightarrow} \check{X}_K(\mathbb{C}). \leqno{(2.1.12)}$$
We call $\check{X}_K(\mathbb{C})$ the {\it Deninger space} associated to $X_K = {\rm Spec}({\cal O}_K)$.\\
\\
{\bf Remark 2.1.13.}  In [D7; Section 4], Deninger requires some conditions on $P^{\times}$ in (2.1.5) to make $\dot{X}_K(\mathbb{C}) = \dot{X}_{\overline{\mathbb{Q}}}(\mathbb{C})/{\rm Gal}(\overline{\mathbb{Q}}/K)$ and hence $\check{X}_K(\mathbb{C})$ smaller space which would lead to a correct geometric model for $X_K$ (cf. [D7; Section 4]). We omit this refinement on $\dot{X}_K(\mathbb{C})$ and $\check{X}_K(\mathbb{C})$, since our purpose in this paper is to give a relation  with Connes-Consani's adelic spaces. \\
\\
{\bf 2.2  Foliated dynamical systems and closed orbits.} We keep the same notations as in 2.1.  Recall that the Deninger space $\check{X}_K(\mathbb{C})$ is equipped with the right Frobenius-$\mathbb{Q}_+$-action as in (2.1.10). We  then define
$$ \frak{X}_K := \check{X}_K(\mathbb{C}) \times_{\mathbb{Q}_+} \mathbb{R}_+.  \leqno{(2.2.1)}$$
to be the quotient of  $\check{X}_K(\mathbb{C}) \times \mathbb{R}_+$ by the {\it right} $\mathbb{Q}_+$-action defined as follows (cf. Notations):
$$ q.((\frak{p},P), u) :=  (F_{q}((\frak{p}, P)), q^{-1}u) = ((\frak{p}, P \circ F_{q}) , q^{-1} u) \; \mbox{for} \; q \in \mathbb{Q}_+.  $$
We equip $\frak{X}_K$ twith he quotient topology of the product $\check{X}_K(\mathbb{C}) \times \mathbb{R}_+$. 

We denote by $[(\frak{p},P), u]$ the image of  $((\frak{p},P), u)$ in $\frak{X}_K$. We then define the Frobenius $\mathbb{R}_+$-flow $\phi$ on $\frak{X}_K$ by the suspension 
$$ \phi^t([(\frak{p}, P), u]) := [(\frak{p},P), t u]. \leqno{(2.2.2)}$$
The 1-codimensional foliation structure on $\frak{X}_K$ is given by taking the images of $\check{X}_K(\mathbb{C})  \times \{ u \}$ $(u \in \mathbb{R}_+)$ as leaves. We call the foliated dynamical system $\frak{X}_K$ the {\it Deninger  system} associated to $X_K = {\rm Spec}({\cal O}_K)$. (In this note, we ignore some refinement on $\check{X}_K(\mathbb{C})$ as explained in Remark 2.1.13.) 

For a Galois extension $L$ of  $K$ in $\overline{\mathbb{Q}}$, the Galois group ${\rm Gal}(L/K)$ acts on $\frak{X}_L$ from the right by
$$ [(\frak{p}, P), u]^{\sigma} := [ (\frak{p}, P)^{\sigma}, u] \;\; \mbox{for} \; \sigma \in {\rm Gal}(L/K).  $$
It commutes with the $\mathbb{R}_+$-action in (2.2.2). The map $\check{\pi}_{L/K}$ in (2.1.11) induces an $\mathbb{R}_+$-equivariant continuous map
$$ \varpi_{L/K} : \frak{X}_L \longrightarrow \frak{X}_K.$$
By (2.1.12), it induces an $\mathbb{R}_+$-equivariant homeomorphism
$$  \frak{X}_L/{\rm Gal}(L/K) \stackrel{\sim}{\longrightarrow}  \frak{X}_K. $$

Next we describe the structure of $\mathbb{R}_+$-orbits in $\frak{X}_K$. By the definition (2.2.1) of $\frak{X}_K$,  it suffices to describe $\mathbb{Q}_+$-orbits in $\check{X}_K(\mathbb{C})$.  Let $\dot{{\rm pr}} : \dot{X}_K(\mathbb{C}) \rightarrow X_K$ and $\check{{\rm pr}} : \check{X}_K(\mathbb{C}) \rightarrow X_K$ be the  projections defined by  $(\frak{p}, P) \rightarrow \frak{p}$, and let $\frak{pr}_K : \frak{X}_K \rightarrow X_K$ be the composition of the projection $\frak{X}_K \rightarrow \check{X}_K(\mathbb{C})$ with $\check{\rm pr}_K$.  These projections $\dot{\rm pr}_K, \check{\rm pr}_K$ and $\frak{pr}_K$ are continuous. By (2.1.10), the $\mathbb{Q}_+$-orbit containing $(\frak{p}, P) \in \check{X}_K(\mathbb{C})$ is mapped to $\frak{p} \in X_K$ under $\check{\rm pr}_K$. So we analyze the structure of the fiber $\check{\rm pr}_K ^{-1}(\frak{p})$. 

We fix a closed point $\frak{p} \in X_K$ with residue characteristic $p$. We take a Galois extension $L$ of $K$, which satisfies the condition  that $$ \begin{array}{ll} \mbox{ any point} \; \frak{Q} \in X_L\; \mbox{over}\; \frak{p} \; \mbox{ has algebraically closed residue } \\
 \; \mbox{field} \; \kappa(\frak{Q}) =   \overline{\mathbb{F}}_p. \end{array} \leqno{(2.2.3)} $$
Since $\overline{\mathbb{F}}_p = \mathbb{F}_p(\mu_{(p)})$,  the smallest Galois extension $L$ of $K$ satisfying the condition (2.2.3) is $L = K(\mu_{(p)})$. Consider the following commutative diagram
\begin{center}
\begin{tikzpicture}[>=stealth,  hookarrow/.style={{Hooks[right]}->},]
\node (0) at (-4, 1) {(2.2.4)};
\node (1) at (-1.2, 2)  {$X_L$};
\node (2) at (-1.2, 0)  {$X_K$};
\node (3) at (2.2, 2)  {$\dot{X}_L(\mathbb{C})$};
\node (4) at (2.2,  0)  {$\dot{X}_K(\mathbb{C})$};
\node (5) at (4.6, 2)  {$\check{X}_L(\mathbb{C})$};
\node (6) at (4.6, 0)  {$\check{X}_K(\mathbb{C})$};
\node (7) at (6.8, 2)  {$\frak{X}_L$} ;
\node (8) at (6.8, 0)  {$\frak{X}_K$} ;
 \draw [->] (1)--(2)node[midway, left]{$\pi_{L/K}$};
  \draw [->] (3)--(4)node[midway, left]{$\dot{\pi}_{L/K}$};
  \draw [->] (5)--(6)node[midway, right]{$\check{\pi}_{L/K}$};
  \draw [->] (7)--(8)node[midway, right]{$\varpi_{L/K}$};
   \draw [->] (3)--(1)node[midway, above]{$\dot{\rm pr}_L$};
   \draw [->] (4)--(2)node[midway,above]{$\dot{\rm pr}_K$};
  \draw [hookarrow] (3)--(5)node[midway, above]{$\;$};
   \draw [hookarrow] (4)--(6)node[midway, above]{$\;$};
   \draw [hookarrow] (5)--(7)node[midway, above]{$\;$};
   \draw [hookarrow] (6)--(8)node[midway, above]{$\;$};
  \end{tikzpicture}
\end{center}
We set
$$\frak{C}_{\frak{p}} := \check{\rm pr}_K ^{-1}(\frak{p}) = \bigcup_{n \in \mathbb{N}} F_n^{-1}(\dot{\rm pr}_K^{-1}(\frak{p})),$$
which is a $\mathbb{Q}_+$-invariant subset of $\check{X}_K(\mathbb{C})$, and we set
$$ \Gamma_{\frak{p}} := \frak{C}_{\frak{p}} \times_{\mathbb{Q}_+} \mathbb{R}_+   \;  \subset \; \frak{X}_K.$$
We note that $\Gamma_{\frak{p}}$ (resp. $\frak{C}_{\frak{p}}$) is the set of $\mathbb{R}_+$-orbits (resp. $\mathbb{Q}_+$-orbits)
passing through a point lying over $\frak{p} \in X_K$ under  $\frak{pr}_K$ (resp. $\check{\rm pr}_K$). We call   $\Gamma_{\frak{p}}$ (resp.  $\frak{C}_{\frak{p}}$) the {\it packet of $\mathbb{R}_+$-orbits} (resp. the {\it packet of $\mathbb{Q}_+$-orbits}) {\it over $\frak{p}$}. By the diagram  (2.2.4), we find the decomposition of $\Gamma_{\frak{p}}$ in $\frak{X}_L$
$$ \varpi_{L/K}^{-1}(\Gamma_{\frak{p}}) = \bigsqcup_{\frak{P}\mid \frak{p}} \Gamma_{\frak{P}}, \leqno{(2.2.5)}$$
where $\frak{P}$ ranges over primes of $L$ over $\frak{p}$. 

Next we describe the structure of each packet  $\frak{C}_{\frak{p}}$ and $\Gamma_{\frak{p}}$. First, using the diagram (2.2.4),  the structure of the  $\mathbb{N}$-invariant set $\dot{\rm pr}_K^{-1}(\frak{p})$ is described as follows. 
By Proposition 2.1.9, we have
$$ \dot{\rm pr}^{-1}(\frak{p}) = \{ (\frak{Q}, P^{\times}) \in \dot{X}_L(\mathbb{C}) \; | \; \pi_{L/K}(\frak{Q}) = \frak{p} \}/{\rm Gal}(L/K).$$ 
Note that ${\rm Gal}(L/K)$ acts transitively on the set of points  $\frak{Q} \in X_L$ over $\frak{p}$. Fix a point $\frak{P} \in X_L$ over $\frak{p}$ so that $\kappa(\frak{P}) = \overline{\mathbb{F}}_p$. Let ${\rm Gal}(L/K)_{\frak{P}}$ be the decomposition group of $\frak{P}$ in ${\rm Gal}(L/K)$, which surjects onto ${\rm Gal}(\kappa(\frak{P})/\kappa(\frak{p})) = {\rm Gal}(\overline{\mathbb{F}}_p/\kappa(\frak{p}))$. 
 Then we have canonical Frobenius $\mathbb{N}$-equivariant bijections
$$  \begin{array}{ll} \dot{\rm pr}^{-1}(\frak{p})   & \simeq \{ (\frak{P}, \overline{P}^{\times}) \in \dot{X}_L(\mathbb{C}) \}/{\rm Gal}(L/K)_{\frak{P}} \\
                                 & \simeq {\rm Hom}_{\bf Gr}(\kappa(\frak{P})^{\times}, \mathbb{C}^{\times})/ {\rm Gal}(\kappa(\frak{P})/\kappa(\frak{p})).
  \end{array}      \leqno{( 2.2.6)}                        
$$
Note that there is the isomorphism
$$ \hat{\mathbb{Z}}_{(p)}^{\times} \stackrel{\sim}{\longrightarrow} {\rm Aut}(\overline{\mathbb{F}}_p^{\times}); \;\; a \mapsto (z \mapsto z^a)$$
and hence there is the natural inclusion
$$ {\rm lk}_{\frak{p}} :    {\rm Gal}(\kappa(\frak{P})/\kappa(\frak{p})) \hookrightarrow {\rm Aut}(\kappa(\frak{P})^{\times}) = \hat{\mathbb{Z}}_{(p)}^{\times}. $$
We note that when $K = \mathbb{Q}$ and $\frak{p} = (p)$, the map ${\rm lk}_{\frak{p}}$  is nothing but the arithmetic linking homomorphism ${\rm lk}_p$ in (1.2.1), (1.2.2) or (1.2.3). More generally, when $K$ is a finite extension of $\mathbb{Q}$, we have ${\rm lk}_{\frak{p}} = {\rm lk}_p^{{\rm deg}(\frak{p})}$ with ${\rm deg}(\frak{p}) := [\kappa(\frak{p}) : \mathbb{F}_p]$. 

We let the monoid $\hat{\mathbb{Z}}_{(p)}^{\times} \times \mathbb{N}$ act on  ${\rm Hom}_{\bf Gr}(\kappa(\frak{P})^{\times}, \mathbb{C}^{\times})$ by pre-composition of $ (\; )^{a} \circ (\; )^n$ $((a, n) \in \hat{\mathbb{Z}}_{(p)}^{\times} \times \mathbb{N})$. 
 The reduction map ${\cal O}_{X_L, \frak{P}} \rightarrow \kappa(\frak{P})$ induces the isomorphism $i_{\frak{P}} : \mu_{(p)}(L) = \mu_{(p)}({\cal O}_{X_L,\frak{P}}) \stackrel{\sim}{\rightarrow} \kappa(\frak{P})^{\times}$. We define $\chi_{\frak{P}}$ by the composition $i_{\frak{P}}^{-1}$ with the inclusion 
$\mu_{(p)}(L)  \hookrightarrow \mathbb{C}^{\times}$:
$$ \chi_{\frak{P}} \; : \; \kappa(\frak{P})^{\times} \stackrel{i_{\frak{P}}^{-1}}{\longrightarrow}    \mu_{(p)}(L) \hookrightarrow \mathbb{C}^{\times}.$$
Then we have a $\hat{\mathbb{Z}}_{(p)}^{\times} \times \mathbb{N}$-equivariant surjection
$$ \hat{\mathbb{Z}}_{(p)}^{\times} \times \mathbb{N}   \twoheadrightarrow {\rm Hom}_{\bf Gr}(\kappa(\frak{P})^{\times}, \mathbb{C}^{\times}); \; (a, n) \mapsto \chi_{\frak{P}}\circ (\; )^{a} \circ (\; )^n.  \leqno{(2.2.7)}$$
Here two elements $(a,n)$ and $(a',n')$ are in the same fiber of this map if and only if $n' = np^m, a = p^m a'$ for some $m \in \mathbb{Z}$. \\
Since ${\rm Gal}(\kappa(\frak{P})/\kappa(\frak{p}))$ is a closed subgroup of $\hat{\mathbb{Z}} \simeq {\rm Gal}(\overline{\mathbb{F}}_p/\mathbb{F}_p)$, we have the following two cases.\\
(i) The case that $\kappa(\frak{p}) = \kappa(\frak{P}) = \overline{\mathbb{F}}_p$. For example, this is the case when $K=L$. By (2.2.6) and (2.2.7), we have the $\mathbb{N}$-equivariant surjection
 $$ \hat{\mathbb{Z}}_{(p)}^{\times} \times \mathbb{N}  \twoheadrightarrow \dot{\rm pr}^{-1}(\frak{p});
 \;\; (a, n) \mapsto \chi_{\frak{P}}\circ (\; )^{a} \circ (\; )^n, $$
where two elements $(a,n)$ and $(a',n')$ are in the same fiber of this map if and only if $n' = np^m, a = p^m a'$ for some $m \in \mathbb{Z}$. Passing to $\varinjlim_{n \in \mathbb{N}}$, we have the following description for the packets $\frak{C}_{\frak{p}}$ and $\Gamma_{\frak{p}}$.   \\
\\
{\bf Theorem 2.2.8} ([D7; Sections 5, 6, 7]). (1) {\it We have the $\mathbb{Q}_+$-equivariant homeomorphism
$$ \hat{\mathbb{Z}}_{(p)}^{\times} \times_{   p^{\mathbb{Z}}} \mathbb{Q}_+ \stackrel{\sim}{\longrightarrow} \frak{C}_{\frak{p}}.$$
and the $\mathbb{R}_+$-equivariant homeomorphism
$$ \hat{\mathbb{Z}}_{(p)}^{\times} \times_{ p^{\mathbb{Z}}} \mathbb{R}_+    \stackrel{\sim}{\longrightarrow} \Gamma_{\frak{p}},$$
 which is the mapping torus of the multiplication by ${\rm lk}_p(p)$ in $\hat{\mathbb{Z}}_{(p)}^{\times}.$ In particular, $\frak{C}_{\frak{p}}$ and $\Gamma_{\frak{p}}$ are independent of $\frak{p}$ over $p$.}\\
 (2)   {\it  Let $\gamma_{\frak{p}, \overline{a}}$ be the fiber in $\Gamma_{\frak{p}}$ over $\overline{a} \in \hat{\mathbb{Z}}_{(p)}^{\times}/p^{\hat{\mathbb{Z}}}$  under the composite of the projections $\Gamma_{\frak{p}} \rightarrow  \hat{\mathbb{Z}}_{(p)}^{\times}/{\rm N}\frak{p}^{\hat{\mathbb{Z}}} \rightarrow \hat{\mathbb{Z}}_{(p)}^{\times}/p^{\hat{\mathbb{Z}}}$.    Then we have a homeomorphism}
 $$ \gamma_{\frak{p}, \overline{a}} \stackrel{\sim}{\longrightarrow} \mathbb{R}_+; \; \; [a, r] \mapsto r,$$
{\it and we have the decomposition into connected $\mathbb{R}_+$-orbits}
$$\Gamma_{\frak{p}} = \bigsqcup_{a \in \hat{\mathbb{Z}}_{(p)}^{\times}} \gamma_{\frak{p},a}.$$ 
\\
(ii) The case that $\kappa(\frak{p})$ is a finite extension of $\mathbb{F}_p$. Equivalently, it is the case that $K$ is a finite extension of $\mathbb{Q}$. Let $f := [\kappa(\frak{p}):\mathbb{F}_p]$. Then we have
$${\rm Gal}(\kappa(\frak{P})/\kappa(\frak{p})) \simeq ({\rm N}\frak{p})^{\hat{\mathbb{Z}}} = p^{f\hat{\mathbb{Z}}} \subset \hat{\mathbb{Z}}_{(p)}^{\times}.$$
 By (2.2.6) and (2.2.7), we have the $\mathbb{N}$-equivariant surjection
 $$ (\hat{\mathbb{Z}}_{(p)}^{\times}/{\rm N}\frak{p}^{\hat{\mathbb{Z}}}) \times \mathbb{N}  \twoheadrightarrow \dot{\rm pr}^{-1}(\frak{p});
 \;\; (\overline{a}, n) \mapsto \chi_{\frak{P}}\circ (\; )^{a} \circ (\; )^n, $$
where $\overline{a} = a {\rm N}\frak{p}^{\hat{\mathbb{Z}}}$ and two elements $(\overline{a},n)$ and $(\overline{a'},n')$ are in the same fiber of this map if and only if $n' = np^m, \overline{a} = p^m \overline{a'}$ for some $m \in \mathbb{Z}$. The latter  depends only on the class mod $f$. Passing to $\varinjlim_{n \in \mathbb{N}}$, we obtain the following  description for $\frak{C}_{\frak{p}}$ and $\Gamma_{\frak{p}}$. \\
\\
{\bf Theorem 2.2.9} ([D7, Sections 5, 6, 7]). (1) {\it We have the $\mathbb{Q}_+$-equivariant homeomorphism
$$ \hat{\mathbb{Z}}_{(p)}^{\times}/{\rm N}\frak{p}^{\hat{\mathbb{Z}}} \times_{p^{\mathbb{Z}}} \mathbb{Q}_+   \stackrel{\sim}{\longrightarrow} \frak{C}_{\frak{p}},$$
and  the $\mathbb{R}_+$-equivariant homeomorphism
$$ \hat{\mathbb{Z}}_{(p)}^{\times}/{\rm N}\frak{p}^{\hat{\mathbb{Z}}}  \times_{p^{\mathbb{Z}}}  \mathbb{R}_+   \stackrel{\sim}{\longrightarrow} \Gamma_{\frak{p}},$$
which is  the mapping torus of the multiplication by ${\rm lk}_p(p)$ in $\hat{\mathbb{Z}}_{(p)}^{\times}/{\rm N}\frak{p}^{\hat{\mathbb{Z}}}$. 
Since ${\rm N}\frak{p} = p^f$ is same for any $\frak{p}$ over $p$, $\frak{C}_{\frak{p}}$ and $\Gamma_{\frak{p}}$ are independent of $\frak{p}$ over $p$.}\\
(2) {\it  Let $\gamma_{\frak{p}, \overline{a}}$ be the fiber in $\Gamma_{\frak{p}}$ over $\overline{a} \in \hat{\mathbb{Z}}_{(p)}^{\times}/p^{\hat{\mathbb{Z}}}$ under the composite of the projections $\Gamma_{\frak{p}} \rightarrow  \hat{\mathbb{Z}}_{(p)}^{\times}/{\rm N}\frak{p}^{\hat{\mathbb{Z}}} \rightarrow \hat{\mathbb{Z}}_{(p)}^{\times}/p^{\hat{\mathbb{Z}}}$.  Then we have the homeomorphism}
 $$ \gamma_{\frak{p}, \overline{a} } \stackrel{\sim}{\longrightarrow} \mathbb{R}_+/{\rm N}\frak{p}^{\mathbb{Z}}; \; \; [a{\rm N}\frak{p}^{\hat{\mathbb{Z}}} , r] \mapsto \overline{r} = r {\rm N}\frak{p}^{\mathbb{Z}} .$$
{\it So $\gamma_{\frak{p}, \overline{a}}$ is a circle of length $\log {\rm N}\frak{p}$ and we have the decomposition into connected closed $\mathbb{R}_+$-orbits}
$$\Gamma_{\frak{p}} = \bigsqcup_{\overline{a} \in \hat{\mathbb{Z}}_{(p)}^{\times}/p^{\hat{\mathbb{Z}}} } \gamma_{\frak{p}, \overline{a} }.$$ 
\\
\\
{\bf 2.3  Monodromies in abelian coverings of Deninger systems.} In this subsection, we consider monodromies in coverings of Deninger systems associated to abelian extensions of $\mathbb{Q}$. To be precise, 
fix a prime number $p$ and let $L$ be a Galois extension of $\mathbb{Q}$, which contains $\mathbb{Q}(\mu_{(p)})$, and let $F$ be a finite subextension of $\mathbb{Q}(\mu_{(p)})$ over $\mathbb{Q}$.  Let $\Gamma_p$ be the packet  of closed orbits in $\frak{X}_{\mathbb{Q}}$ over $p$. We then have the following coverings of Deninger systems and inverse images of $\Gamma_p$:
\begin{center}
\begin{tikzpicture}[>=stealth,  hookarrow/.style={{Hooks[right]}->},]
\node (1) at (4, 3.8)  {$\frak{X}_{L}$} ;
\node (2) at (4, 1.6)  {$\frak{X}_F$} ;
\node (3) at (4, 0.2)  {$\frak{X}_{\mathbb{Q}}$} ;
\node (4) at (4.8, 3.8)  {$\supset$} ;
\node (5) at (6.2, 3.8)  {$\varpi_{L/\mathbb{Q}}^{-1}(\Gamma_p)$}; 
\node (7) at (8.2, 3.8) {$= \bigsqcup_{\frak{P}\mid p}\Gamma_{\frak{P}}$} ;
\node (8) at (6.2, 1.6)  {$\varpi_{F/\mathbb{Q}}^{-1}(\Gamma_p)$};
\node (9) at (8, 1.6) {$= \bigsqcup_{\frak{p}\mid p}\Gamma_{\frak{p}}$} ;
\node (12) at (5, 1.6)  {$\supset$} ;
\node (11) at (6.2, 0.2)  {$\Gamma_{p}$} ;
\node (10) at (5, 0.2)  {$\supset$} ;
   \draw [->] (1)--(2)node[midway, left]{$\varpi_{L/F}$};
  \draw [->] (2)--(3)node[midway,  left]{$\varpi_{F/\mathbb{Q}}$};
   \draw [->] (5)--(8)node[midway, right]{$\;$};
  \draw [->] (8)--(11)node[midway,  right]{$\;$};
\end{tikzpicture}
\end{center}
Here $\Gamma_{\frak{P}}$ (resp. $\Gamma_{\frak{p}}$) is a packet of $\mathbb{R}_+$-orbits and $\frak{P}$ (resp. $\frak{p}$) ranges over primes of $L$ (resp. $F$) over $p$ (cf. (2.2.3)). By Theorem 2.2.9 (2), the packet  $\Gamma_p$ consists of circles $\gamma_{p,  \overline{a}} = \mathbb{R}_+/p^{\mathbb{Z}}$, which are fibers over $\overline{a} \in \hat{\mathbb{Z}}_{(p)}^{\times}/p^{\hat{\mathbb{Z}}}$. In this subsection, we describe the monodromies around $\gamma_{p, \overline{a}}$ in $\varpi_{L/\mathbb{Q}}^{-1}(\gamma_{p, \overline{a}})$ and $\varpi_{F/\mathbb{Q}}^{-1}(\gamma_{p, \overline{a}})$, according as the cases (i) and (ii) in the subsection 2.2..  \\
\\
(i)  By Theorem 2.2.8 (1), we have the homeomorphism
$$ \Gamma_{\frak{P}} \simeq \hat{\mathbb{Z}}_{(p)}^{\times} \times_{p^{\mathbb{Z}}} \mathbb{R}_+.$$
Looking at the commutative  diagram
\begin{center}
\begin{tikzpicture}[>=stealth,  hookarrow/.style={{Hooks[right]}->},]
\node (1) at (-1.2, 2)  {$\Gamma_{\frak{P}}$};
\node (2) at (-1.2, 0)  {$\Gamma_p$};
\node (3) at (2, 2)  {$\hat{\mathbb{Z}}_{(p)}^{\times}$};
\node (4) at (2,  0)  {$\hat{\mathbb{Z}}_{(p)}^{\times}/p^{\hat{\mathbb{Z}}},$};
 \draw [->] (1)--(2)node[midway, left]{$\varpi_{L/\mathbb{Q}}|_{\Gamma_{\frak{P}}}$};
  \draw [->] (3)--(4)node[midway, right]{mod $p^{\hat{\mathbb{Z}}}$};
    \draw [->] (1)--(3)node[midway, above]{$\;$};
   \draw [->] (2)--(4)node[midway,above]{$\;$};
  \end{tikzpicture}
\end{center}
we see that $(\varpi_{L/\mathbb{Q}}|_{\Gamma_{\frak{P}}})^{-1}(\gamma_{p, \overline{a}})$ is the fiber $\gamma_{\frak{P}, \overline{a}}$ in $\Gamma_{\frak{P}}$ over $\overline{a} \in \hat{\mathbb{Z}}_{(p)}^{\times}/p^{\hat{\mathbb{Z}}}$ under the composite of the projections
 $\Gamma_{\frak{P}} \rightarrow \hat{\mathbb{Z}}_{(p)}^{\times} \rightarrow \hat{\mathbb{Z}}_{(p)}^{\times}/p^{\hat{\mathbb{Z}}}$ (cf. Theorem 2.2.9 (2)). Here it may be noted  that the inverse image of $\overline{a}$ in $\hat{\mathbb{Z}}_{(p)}^{\times}$ is $\{ ap^z \mid z \in \hat{\mathbb{Z}}\}$ and the fibers in $\Gamma_{\frak{P}}$ over $ap^z$  are  same 
for $z  \in \hat{\mathbb{Z}}$, which is $\gamma_{\frak{p}, \overline{a}}$.
Therefore we have 
$$ \varpi_{L/\mathbb{Q}}^{-1}(\gamma_{p, \overline{a}}) = \bigsqcup_{\frak{P} \mid p} \gamma_{\frak{P}, \overline{a}}.$$
By Theorem 2.2.8 (2) and Theorem 2.2.9 (2), we have the commutative diagram
\begin{center}
\begin{tikzpicture}[>=stealth,  hookarrow/.style={{Hooks[right]}->},]
\node (1) at (-1.2, 2)  {$\gamma_{\frak{P}, \overline{a}}$};
\node (2) at (-1.2, 0)  {$\gamma_{p, \overline{a}}$};
\node (3) at (2, 2)  {$\mathbb{R}_+$};
\node (4) at (2,  0)  {$\mathbb{R}_+/p^{\mathbb{Z}}$};
 \draw [->] (1)--(2)node[midway, left]{$\varpi_{L/\mathbb{Q}}|_{\gamma_{\frak{P}, \overline{a}}}$};
  \draw [->] (3)--(4)node[midway, right]{mod $p^{\mathbb{Z}}$};
    \draw [->] (1)--(3)node[midway, above]{$\sim$};
   \draw [->] (2)--(4)node[midway,above]{$\sim$};
  \end{tikzpicture}
\end{center}
Here the monodromy around $\gamma_{p, \overline{a}}$ in $\gamma_{\frak{P}, \overline{a}}$ is the multiplication by $p$ in 
$\mathbb{R}_+$.  Summarizing the above, we have the following\\
\\
{\bf Theorem 2.3.1.} {\it Notations being as above,   the inverse image $\varpi_{L/\mathbb{Q}}^{-1}(\gamma_{p, \overline{a}})$ of $\gamma_{p, \overline{a}}$ is decomposed  into the connected components of $\mathbb{R}_+$-orbits}
$$ \varpi_{L/\mathbb{Q}}^{-1}(\gamma_{p, \overline{a}}) = \bigsqcup_{\frak{P} \mid p} \gamma_{\frak{P}, \overline{a}}, \;\; \gamma_{\frak{P}, \overline{a}} \simeq \mathbb{R}_+.$$
{\it The monodromy around $\gamma_{p, \bar{a}}$ in the fiber $\varpi_{L/\mathbb{Q}}^{-1}(\gamma_{p, \bar{a}})$ is given by the multiplication by $p$ in each connected component $(\mathbb{R}_+$-orbit$)$.}\\
\\
(ii) By Theorem 2.2.9 (1), for each prime $\frak{p}$ of $F$ over $p$, we have the homeomorphism
$$ \Gamma_{\frak{p}} \simeq \hat{\mathbb{Z}}_{(p)}^{\times}/{\rm N}\frak{p}^{\hat{\mathbb{Z}}} \times_{p^{\mathbb{Z}}} \mathbb{R}_+.$$
Let $f := [\kappa(\frak{p}) : \mathbb{F}_p]$. Since $p$ is unramified in $F$, the cyclic group $p^{\hat{\mathbb{Z}}}/{\rm N}\frak{p}^{\hat{\mathbb{Z}}} \simeq \mathbb{Z}/f \mathbb{Z}$ is the decomposition group of $\frak{p}$ in $F/\mathbb{Q}$. 
Looking at the commutative  diagram
\begin{center}
\begin{tikzpicture}[>=stealth,  hookarrow/.style={{Hooks[right]}->},]
\node (1) at (-1, 2)  {$\Gamma_{\frak{p}}$};
\node (2) at (-1, 0)  {$\Gamma_p$};
\node (3) at (2, 2)  {$\hat{\mathbb{Z}}_{(p)}^{\times}/{\rm N}\frak{p}^{\hat{\mathbb{Z}}}$};
\node (4) at (2,  0)  {$\hat{\mathbb{Z}}_{(p)}^{\times}/p^{\hat{\mathbb{Z}}},$};
 \draw [->] (1)--(2)node[midway, left]{$\varpi_{F/\mathbb{Q}}|_{\Gamma_{\frak{p}}}$};
  \draw [->] (3)--(4)node[midway, right]{mod $p^{\hat{\mathbb{Z}}}$};
    \draw [->] (1)--(3)node[midway, above]{$\;$};
   \draw [->] (2)--(4)node[midway,above]{$\;$};
  \end{tikzpicture}
\end{center}
we see that $(\varpi_{F/\mathbb{Q}}|_{\Gamma_{\frak{p}}})^{-1}(\gamma_{p, \overline{a}})$ is the fiber $\gamma_{\frak{p}, \overline{a}}$ in $\Gamma_{\frak{p}}$ over $\overline{a} \in \hat{\mathbb{Z}}_{(p)}^{\times}/p^{\hat{\mathbb{Z}}}$ under the composite of the projections 
$\Gamma_{\frak{p}} \rightarrow \hat{\mathbb{Z}}_{(p)}^{\times}/{\rm N}\frak{p}^{\hat{\mathbb{Z}}} \rightarrow  \hat{\mathbb{Z}}_{(p)}^{\times}/p^{\hat{\mathbb{Z}}} $ (cf. Theorem 2.2.9 (2)). 
It may be noted that the inverse image of $\overline{a}$ in $\hat{\mathbb{Z}}_{(p)}^{\times}/{\rm N}\frak{p}^{\hat{\mathbb{Z}}}$ is $\{ ap^z {\rm N}\frak{p}^{\hat{\mathbb{Z}}} \mid 
z \in  \mathbb{Z}/f\mathbb{Z} \}$ and the fibers in  $\Gamma_{\frak{p}}$ over $a p^z {\rm N}\frak{p}^{\hat{\mathbb{Z}}}$ are same for $z \in \mathbb{Z}/f\mathbb{Z}$, which is $\gamma_{\frak{p}, \overline{a}}$.
By Theorem 2.2.9 (2), we have the commutative diagram
\begin{center}
\begin{tikzpicture}[>=stealth,  hookarrow/.style={{Hooks[right]}->},]
\node (1) at (-1.2, 2)  {$ \gamma_{\frak{p}, \overline{a}}$};
\node (2) at (-1.2, 0)  {$\gamma_{p, \overline{a}}$};
\node (3) at (2, 2)  {$\mathbb{R}_+/{\rm N}\frak{p}^{\mathbb{Z}}$};
\node (4) at (2,  0)  {$\mathbb{R}_+/p^{\mathbb{Z}}$};
 \draw [->] (1)--(2)node[midway, left]{$\varpi_{F/\mathbb{Q}}|_{\gamma_{\frak{p}, \overline{a p^z}}}$};
  \draw [->] (3)--(4)node[midway, right]{mod $p^{\mathbb{Z}}$};
    \draw [->] (1)--(3)node[midway, above]{$\sim$};
   \draw [->] (2)--(4)node[midway,above]{$\sim$};
  \end{tikzpicture}
\end{center}
Here the monodromy around $\gamma_{p, \overline{a}}$ in $\gamma_{\frak{p}, \overline{a}}$ is the multiplication by $p$ in 
$\mathbb{R}_+/{\rm N}\frak{p}^{\mathbb{Z}}$.  Now we have
$$  \varpi_{F/\mathbb{Q}}^{-1}(\gamma_{p, \overline{a}}) = \bigsqcup_{\frak{p} \mid p} \gamma_{\frak{p}, \overline{a}} = \bigsqcup_{\sigma \in {\rm Gal}(F/\mathbb{Q})} \gamma_{\frak{p}^{\sigma}, \overline{a}}, $$
where we fix a prime $\frak{p}$ over $p$ on the right hand side. Let
$$ \chi_F : \hat{\mathbb{Z}}_{(p)}^{\times} = {\rm Gal}(\mathbb{Q}(\mu_{(p)})/\mathbb{Q}) \longrightarrow {\rm Gal}(F/\mathbb{Q})
$$
be the restriction homomorphism so that 
$$ \chi_F({\rm lk}_p(p)) = \chi_F(\sigma_p) = \left( \frac{F/\mathbb{Q}}{p} \right) $$
is the Artin symbol over $p$ in $F/\mathbb{Q}$ as in (1.2.5). Then we have the bijection
$$ \bigsqcup_{\sigma \in {\rm Gal}(F/\mathbb{Q})} \gamma_{\frak{p}^{\sigma}, \overline{a}}
\stackrel{\sim}{\longrightarrow} {\rm Gal}(F/\mathbb{Q}) \times_{p^{\mathbb{Z}}}  \mathbb{R}_+; \; \; [a {\rm N}\frak{p}^{\sigma}, r] \mapsto [\chi_F(a), r].$$
So the fiber $\varpi_{F/\mathbb{Q}}^{-1}(\gamma_{p, \overline{a}})$ is described as the mapping torus ${\rm Gal}(F/\mathbb{Q}) \times_{p^{\mathbb{Z}}}  \mathbb{R}_+$ and the monodromy around the circle $\gamma_{p, \overline{a}}$ is the multiplication by $\displaystyle{\left( \frac{F/\mathbb{Q}}{p} \right)}$ in ${\rm Gal}(F/\mathbb{Q})$. Summarizing the above, we have the following\\
\\
{\bf Theorem 2.3.2.} {\it Notations being as above, the following assertions hold.} \\
(1) {\it  The inverse image $\varpi_{F/\mathbb{Q}}^{-1}(\gamma_{p, \overline{a}})$ of $\gamma_{p, \overline{a}}$ is decomposed into the connected components of $\mathbb{R}_+$-orbits}\\
$$ \varpi_{F/\mathbb{Q}}^{-1}(\gamma_{p, \overline{a}})  = \bigsqcup_{\frak{p} \mid p} \gamma_{\frak{p}, \overline{a}}, \;\; \gamma_{\frak{p}, \overline{a}} \simeq \mathbb{R}_+/{\rm N}\frak{p}^{\mathbb{Z}}$$
{\it and the monodromy around $\gamma_{p, \overline{a}}$ in $\varpi_{F/\mathbb{Q}}^{-1}(\gamma_{p, \overline{a}})$ is the multiplication by
$p$ in each connected component $(\mathbb{R}_+$-orbit$)$. {\it  We can also describe  $\varpi_{F/\mathbb{Q}}^{-1}(\gamma_{p, \overline{a}})$ as the mapping torus}
$$\varpi_{F/\mathbb{Q}}^{-1}(\gamma_{p, \overline{a}})   \simeq {\rm Gal}(F/\mathbb{Q}) \times_{p^{\mathbb{Z}}}  \mathbb{R}_+, $$
where the monodromy  around $\gamma_{p, \overline{a}}$ in $\varpi_{F/\mathbb{Q}}^{-1}(\gamma_{p, \overline{a}})$  is given by the multiplication by $\displaystyle{\left( \frac{F/\mathbb{Q}}{p} \right)}$ in ${\rm Gal}(F/\mathbb{Q})$.  } \\
(2) {\it Let $(p) = \frak{p}_1 \cdots \frak{p}_r$ be the decomposition of $(p)$ into distinct primes in $F$ such that $f = {\rm deg}(\frak{p}_i) = [\kappa(\frak{p}_i) : \mathbb{F}_p]$ and $fr = [F:\mathbb{Q}]$. Then $\varphi_{F/\mathbb{Q}}^{-1}(\gamma_{p, \overline{a}})$ is decomposed into connected components $($circles$)$ in the same manner, namely, 
$\varpi_{F/\mathbb{Q}}^{-1}(\gamma_{p, \overline{a}}) = \gamma_{\frak{p}_1, \overline{a}}  \sqcup \cdots \sqcup \gamma_{\frak{p}_r, \overline{a}} $
 such that $\varpi_{F/\mathbb{Q}}|_{\gamma_{\frak{p}_i, \overline{a}}} : \gamma_{\frak{p}_i, \overline{a}} = \mathbb{R}_+/{\rm N}\frak{p}^{\mathbb{Z}} \rightarrow \gamma_{p, \overline{a}} = \mathbb{R}_+/p^{\mathbb{Z}}$ is a cyclic covering of degree $f$ and $fr = [F:\mathbb{Q}]$. }\\
\\
Since the assertions in Theorem 2.3.1 and 2.3.2 are independent of $\overline{a} \in \hat{\mathbb{Z}}_{(p)}^{\times}/p^{\hat{\mathbb{Z}}}$,
it makes sense to use the phrase ``the monodromy around $\Gamma_p$ in $\varpi_{K/\mathbb{Q}}^{-1}(\Gamma_p)$" for $K = L, F$.
\\
\\
{\bf Example 2.3.3.} Let $p$ and $q$ be distinct odd prime numbers and let $F = \mathbb{Q}(\sqrt{q})$. Then the monodromy around 
$\Gamma_p$ in  $\varpi_{F/\mathbb{Q}}^{-1}(\Gamma_p)$ is given by the multiplication by the multiplication of $p$ by 
$\bigsqcup_{\frak{p} \mid p} \hat{\mathbb{Z}}_{(p)}^{\times}/{\rm N}\frak{p}^{\hat{\mathbb{Z}}}$.
$$ \begin{array}{ccc}
\varpi_{F/\mathbb{Q}}^{-1}(\Gamma_p) = \Gamma_{\frak{p}} \cup \Gamma_{\overline{\frak{p}}} \; (\mbox{disjoint}) & \Longleftrightarrow & \displaystyle{\left( \frac{q}{p} \right)}= 1,\\
\varpi_{F/\mathbb{Q}}^{-1}(\Gamma_p) = \Gamma_{\frak{p}}   & \Longleftrightarrow & \displaystyle{\left( \frac{q}{p} \right)} = -1.\\
\end{array}
$$
\\
{\bf Remark 2.2.4.}  We note that there are similarities between Theorems 1.2.4, 1.2.6 for Connes-Consani adelic spaces and Theorems 2.3.1, 2.3.2 for Deninger systems in the respect that the arithmetic monodromy over a prime in an extension of number fields is described by the geometric monodromy along a closed orbit in a covering of dynamical systems. \\
\\

\begin{center}
{\bf 3.  A relation between Deninger systems and Connes-Consani adelic spaces}
\end{center}

In this section, we give a relation between Deninger systems for number rings and Connes-Consani adelic spaces. First, we recall Deninger's construction of the maps from Deninger spaces to finite adeles,
and  then we bridge Deninger systems $\frak{K}$ and Connes-Consani adelic spaces $\mathscr{X}_F$ for finite abelian extensions $F$ of $\mathbb{Q}$. We keep the same notations as in the sections 1 and 2. \\

We set
$$\dot{H} := {\rm Hom}_{\bf Gr}(\mu_{\infty}, \mu_{\infty}),$$
which is a topological ring equipped with the topology of pointwise convergence. Then we have an isomorphism of topological rings
$$ \hat{\mathbb{Z}} \stackrel{\sim}{\longrightarrow} \dot{H}; \; a \mapsto (\chi_a : \zeta \mapsto \zeta^a).  \leqno{(3.1)}$$
The $\mathbb{N}$-action on $\dot{H}$  is defined by exponentiation of characters 
$$n. \chi := \chi^n \;\; \;\; (n \in \mathbb{N},  \chi \in \dot{H}). \leqno{(3.2)}$$
It is just the $\mathbb{N}$-action by multiplication on the abelian group $\dot{H}$.
 The Galois group ${\rm Gal}(\mathbb{Q}^{\rm ab}/\mathbb{Q})$ acts continuously on $\dot{H}$ by
$ \chi^{\sigma}(\zeta) := \chi(\zeta^{\sigma})$
for $\chi \in \dot{H}, \sigma \in {\rm Gal}(\mathbb{Q}^{\rm ab}/\mathbb{Q})$ and $\zeta \in \mu_{\infty}$. Let 
$ {\rm cycl} : {\rm Gal}(\mathbb{Q}^{\rm ab}/\mathbb{Q}) \rightarrow \hat{\mathbb{Z}}^{\times} $
be the cyclotomic character defined by $\zeta^{\sigma} = \zeta^{{\rm cycl}(\sigma)} $
for $\zeta \in \mu_{\infty}$ and $\sigma \in {\rm Gal}(\mathbb{Q}^{\rm ab}/\mathbb{Q})$. Then we have
$$ \chi^{\sigma} = \chi^{{\rm cycl}(\sigma)} \;  \mbox{pointwise for} \; \chi \in \dot{H},  \sigma \in {\rm Gal}(\mathbb{Q}^{\rm ab}/\mathbb{Q}). \leqno{(3.3)}$$

Let $(\frak{P}, P) \in \dot{X}_{\mathbb{Q}^{\rm ab}}(\mathbb{C})$. Recall from (2.1.3) that $\frak{P} \in X_{\mathbb{Q}^{\rm ab}}$ and $P : W_{\rm rat}(\kappa(\frak{P})) \rightarrow \mathbb{C}$ is a ring homomorphism. 
If $\frak{P}$ is a closed point of $X_{\mathbb{Q}^{\rm ab}}$ whose residue characteristic is $p$, then $\kappa(\frak{P}) = \overline{\mathbb{F}}_p$. Since $W_{\rm rat}(\kappa(\frak{P})) \simeq \mathbb{Z}[\overline{\mathbb{F}}_p^{\times}]$ by Example 2.1.2, we have the multiplicative inclusion $\mu_{\infty}(\kappa(\frak{P})) = \mu_{(p)} \hookrightarrow W_{\rm rat}(\kappa(\frak{P}))$. If $\frak{P}$ is the generic point of  $X_{\mathbb{Q}^{\rm ab}}$, then $\kappa(\frak{P}) = \mathbb{Q}^{\rm ab}$.  Since we have the ring isomorphism $ \mathbb{Z}[\overline{\mathbb{Q}}^{\times}]^{{\rm Gal}(\overline{\mathbb{Q}}/\mathbb{Q}^{\rm ab})} \simeq W_{\rm rat}(\mathbb{Q}^{\rm ab})$, we also have the multiplicative inclusion $\mu_{\infty}(\kappa(\frak{P})) = \mu_{\infty} \hookrightarrow W_{\rm rat}(\mathbb{Q}^{\rm ab})$. Thus the image of the multiplicative map $P|_{\mu_{\infty}(\kappa(\frak{P}))}$ is in $\mu_{\infty}$ for any $\frak{P} \in X_{\mathbb{Q}^{\rm ab}}$ and any $P \in {\rm Hom}_{\bf CRing}(W_{\rm rat}(\kappa(\frak{P}), \mathbb{C})$.  Let  
$$\iota_{\frak{P}} : \mu_{\infty}  = \mu_{\infty}({\cal O}_{X_{\mathbb{Q}^{\rm ab}}, \frak{P}}) \longrightarrow \mu_{\infty}(\kappa(\frak{P}))$$
be the map obtained by restricting  the reduction map ${\cal O}_{X_{\mathbb{Q}^{\rm ab}}, \frak{P}} \rightarrow \kappa(\frak{P})$
to $\mu_{\infty}({\cal O}_{X_{\mathbb{Q}^{\rm ab}}, \frak{P}} )$ and we set 
$$ \chi_{(\frak{P}, P)} := P|_{\mu_{\infty}(\kappa(\frak{P}))} \circ \iota_{\frak{P}}.$$
Then we have the map
$$ \psi : \dot{X}_{\mathbb{Q}^{\rm ab}}(\mathbb{C}) \longrightarrow \dot{H}; \;\; (\frak{P}, P) \mapsto \chi_{(\frak{P}, P)}.$$
\vspace{0.02cm}\\
{\bf Lemma 3.4.} {\it The map $\psi$ is an $\mathbb{N}$- and ${\rm Gal}(\mathbb{Q}^{\rm ab}/\mathbb{Q})$-equivariant continuous map.}
\begin{proof}
 By (2.1.4), (3.2), (2.1.7) and (3.3), it suffices to show
$$ \psi((\frak{P}, P \circ F_n)) = \chi_{(\frak{P},P)}^n, \;\; \psi((\frak{P}^{\sigma}, P \circ \sigma)) = \chi_{(\frak{P}, P)}^{\rm cycl} $$
for the $\mathbb{N}$- and ${\rm Gal}(\mathbb{Q}^{\rm ab}/\mathbb{Q})$-equivariance of $\psi$.  By evaluating $\zeta \in \mu_{\infty}$, we can easily check these. Since $\psi$ is given by restricting $P$ to $i_{\frak{P}}(\mu_{\infty})$, the continuity of $\psi$ follows. 
\end{proof} 
\par\medskip
\par\noindent
We order $\mathbb{N}$ by division and  define the topological ring $\check{H}$ by
$$ \check{H} := \varinjlim_{n \in \mathbb{N}} \dot{H} = \dot{H} \otimes_{\mathbb{Z}} \mathbb{Q}$$ 
with the inductive limit topology, on which $ {\rm Gal}(\mathbb{Q}^{\rm ab}/\mathbb{Q}) \times \mathbb{Q}_+$ acts. By (3.1), we obtain the isomorphism of topological rings 
$$ \check{H}  \simeq  \hat{\mathbb{Z}} \otimes_{\mathbb{Z}} \mathbb{Q} = \mathbb{A}^{\rm f}.$$
By Lemma 3.4 and passing to the inductive limit over $\mathbb{N}$, we obtain the $\mathbb{Q}_+$- and ${\rm Gal}(\mathbb{Q}^{\rm ab}/\mathbb{Q})$-equivariant continuous map
$$ \check{\psi} : \check{X}_{\mathbb{Q}^{\rm ab}}(\mathbb{C}) \longrightarrow \check{H} \simeq \mathbb{A}^{\rm f} $$
and hence
$$ \check{\psi} \times {\rm id}_{\mathbb{R}_+} : \check{X}_{\mathbb{Q}^{\rm ab}}(\mathbb{C}) \times \mathbb{R}_+ \longrightarrow \mathbb{A}^{\rm f}  \times \mathbb{R}_+ $$
is a $\mathbb{R}_+$- and ${\rm Gal}(\mathbb{Q}^{\rm ab}/\mathbb{Q})$-equivariant continuous map, where $\mathbb{R}_+$ acts on the $\mathbb{R}_+$-component in the both sides by multiplication.  By the definition of the Deninger system, taking the quotient by the left $\mathbb{Q}$-action (2.2.1), we obtain the $\mathbb{R}_+$- and ${\rm Gal}(\mathbb{Q}^{\rm ab}/\mathbb{Q})$-equivariant continuous map
$$ \check{\psi} \times_{\mathbb{Q_+}} {\rm id}_{\mathbb{R}_+} : \frak{X}_{\mathbb{Q}^{\rm ab}} \longrightarrow \mathbb{A}^{\rm f} \times_{\mathbb{Q}_+} \mathbb{R}_+. $$
Let 
$$ \tau : \mathbb{A}^{\rm f} \times_{\mathbb{Q}_+} \mathbb{R}_+ \longrightarrow \mathbb{Q}^{\times} \backslash \mathbb{A} = : \mathscr{X}_{\mathbb{Q}^{\rm ab}}$$
be the composite of the homeomorphism $ \mathbb{A}^{\rm f} \times_{\mathbb{Q}_+} \mathbb{R}_+ \longrightarrow \mathbb{Q}_+ \backslash (\mathbb{A}^{\rm f} \times \mathbb{R}_+)$ defined by $((a_p), r) \mapsto ((a_p), r^{-1})$ with the natural map $\mathbb{Q}_+ \backslash (\mathbb{A}^{\rm f} \times \mathbb{R}_+) \rightarrow \mathbb{Q}^{\times} \backslash \mathbb{A}$. By composing two continuous maps $\psi \times_{\mathbb{Q_+}} {\rm id}_{\mathbb{R}_+}$ and $\tau$, 
we obtain the continuous map
$$ \Psi_{\mathbb{Q}^{\rm ab}} : \frak{X}_{\mathbb{Q}^{\rm ab}}  \longrightarrow \mathscr{X}_{\mathbb{Q}^{\rm ab}}. $$
\\
{\bf Lemma 3.5.} {\it The map $\Psi_{\mathbb{Q}^{\rm ab}}$ is an $\mathbb{R}_+$-anti-equivariant and ${\rm Gal}(\mathbb{Q}^{\rm ab}/\mathbb{Q})$-equivariant continuous map.} 
\begin{proof}  By (1.1.4) and (2.2.2), the $\mathbb{R}_+$-action on $\mathscr{X}_{\mathbb{Q}^{\rm ab}}$ is the multiplication on the 
$\infty$-component and that on $\frak{X}_{\mathbb{Q}^{\rm ab}}$ is the multiplication on the 
$\mathbb{R}_+$-component. So the twisting in $\tau$ makes $\Psi_{\mathbb{Q}^{\rm ab}}$ $\mathbb{R}_+$-anti-equivariant. 
We showed $\check{\psi}$ is ${\rm Gal}(\mathbb{Q}^{\rm ab}/\mathbb{Q})$-equivariant under the ${\rm Gal}(\mathbb{Q}^{\rm ab}/\mathbb{Q})$-action on $\check{H} \otimes \mathbb{Q} \simeq \mathbb{A}^{\rm f}$ defined by (3.3). Since ${\rm Gal}(\mathbb{Q}^{\rm ab}/\mathbb{Q})$ acts on $\mathscr{X}_{\mathbb{Q}^{\rm ab}}$ by (1.1.3), it suffices to prove
the isomorphism $\hat{\mathbb{Z}} \stackrel{\sim}{\rightarrow} \dot{H}; a \mapsto \chi_a $ in (3.1) is compatible with the ${\rm Gal}(\mathbb{Q}^{\rm ab}/\mathbb{Q})$-actions,
where ${\rm Gal}(\mathbb{Q}^{\rm ab}/\mathbb{Q})$ acts on $\dot{H}$ by (3.3) and on $\hat{\mathbb{Z}}$ via the multiplication
of $\hat{\mathbb{Z}}^{\times}$ under the identification ${\rm Gal}(\mathbb{Q}^{\rm ab}/\mathbb{Q}) = \hat{\mathbb{Z}}^{\times}$ by ${\rm cycl}$.
But it follows from that $\chi_a^{\sigma} = \chi_a^{{\rm cycl}(\sigma)} = \chi_{a \cdot{\rm cycl}(\sigma)}$ for $a \in \hat{\mathbb{Z}}, \sigma \in {\rm Gal}(\mathbb{Q}^{\rm ab}/\mathbb{Q})$.  \end{proof} 
\par\medskip
\par\noindent
Let $F$ be a finite subextension of $\mathbb{Q}^{\rm ab}$ over $\mathbb{Q}$. By Lemma 3.5 and taking the quotient by ${\rm Gal}(\mathbb{Q}^{\rm ab}/F) = U_F$ in $\Psi_{\mathbb{Q}^{\rm ab}}$, we obtain the anti-$\mathbb{R}_+$-equivariant and ${\rm Gal}(F/\mathbb{Q}) = {\rm Gal}(\mathbb{Q}^{\rm ab}/F) /U_F$-equivariant continuous map
$$ \Psi_{F} : \frak{X}_{F} \longrightarrow \mathscr{X}_{F}. $$
\\
{\bf Theorem 3.6.} {\it Notations being as above, the following assertions hold.}\\
(1) {\it We have the following commutative diagram, where the horizontal arrows are $\mathbb{R}_+$-anti-equivariant continuous maps.}
\begin{center}
\begin{tikzpicture}[>=stealth]
\node (1) at (0.2, 2.2)  {$\frak{X}_{\mathbb{Q}^{\rm ab}}$} ;
\node (2) at (3.8, 2.2)  {$\mathscr{X}_{\mathbb{Q}^{\rm ab}}$} ;
\node (3) at (0.2, 0)  {$\frak{X}_F$} ;
\node (4) at (3.8, 0)  {$\mathscr{X}_F$} ;
\node (5) at (0.2, -1.6)  {$\frak{X}_{\mathbb{Q}}$} ;
\node (6) at (3.8, -1.6)  {$\mathscr{X}_{\mathbb{Q}}$} ;
 \draw [->] (1)--(2)node[midway,above]{$\Psi_{\mathbb{Q}^{\rm ab}}$};
  \draw [->] (3)--(4)node[midway,above]{$\Psi_F$};
  \draw [->] (1)--(3)node[midway,left]{$\varpi_{\mathbb{Q}^{\rm ab}/F}$};
   \draw [->] (2)--(4)node[midway,right]{$\varphi_{\mathbb{Q}^{\rm ab}/F}$};
   \draw [->] (5)--(6)node[midway,above]{$\Psi_{\mathbb{Q}}$};
  \draw [->] (3)--(5)node[midway,left]{$\varpi_{F/\mathbb{Q}}$};
   \draw [->] (4)--(6)node[midway,right]{$\varphi_{F/\mathbb{Q}}$};   
\end{tikzpicture}
\end{center}
(2) {\it For a prime number $p$, an $\mathbb{R}_+$-orbit $($circle$)$ $\gamma_p$ in the packet $\Gamma_p$ is sent onto the circle $C_p$ under $\Psi_{\mathbb{Q}}$. 
If $p$ is unramified in $F$, $\varpi_{F/\mathbb{Q}}^{-1}(\gamma_p)$ is sent onto $\varphi_{F/\mathbb{Q}}^{-1}(C_p)$.}
\begin{proof} Since the assertions in (1) were already shown, we need to show (2). Let $\gamma_p$ be an $\mathbb{R}_+$-orbit in $\Gamma_p$. Let $\frak{P} \in X_{\mathbb{Q}^{\rm ab}}$ be a prime over $p$ and $(\frak{P}, P) \in \dot{X}_{\mathbb{Q}^{\rm ab}}(\mathbb{C})$. By the definition (1.1.5) of $C_p$ and the commutative diagram in (1), we have only to claim that the $p$-component of $\psi((\frak{P},P)) \in \dot{H} = \hat{\mathbb{Z}}$ is  zero.  Since $\chi_{(\frak{P},P)}(\zeta) = (P\circ \iota_{\frak{P}})(\zeta) = \zeta^a  \in \mu_{(p)}$ for $\zeta \in \mu_{\infty}$, the $p$-component of $a$ is zero, which yields the claim. The last assertion follows from the descriptions of  $\varphi_{F/\mathbb{Q}}^{-1}(C_p)$ in Theorem 1.2.8 and of $\varpi_{F/\mathbb{Q}}^{-1}(\gamma_p)$ in Theorem 2.3.2. \end{proof} 
\par\medskip
\par\noindent
{\bf References}	
\begin{flushleft}
{[BC]} J.-B. Bost and A. Connes, Hecke algebras, type III factors and phase transitions with spontaneous symmetry breaking in number theory, Selecta Mathematica, (N.S.) {\bf 1} (3),  1995, 411--457. \\
{[C1]}  A.  Connes, Noncommutative geometry, Academic Press, 1994.\\
{[C2]} A. Connes, Formule de trace en g\'{e}om\'{e}trie non commutative et hypoth\`{e}ses de Riemann, Comptes Rendus de l'Acad\'{e}mie des Sci.ences, Paris Ser. I.  Math.,, {\bf 323}, no. 12,  1996, 1231--1236.\\
{[C3]} A. Connes, Trace formula in noncommutative geometry and the zeros of the Riemann zeta function, Selecta Mathematica, (N.S) {\bf 5} (1), 1999, 29--106. \\
{[CC1]} A. Connes, C. Consani, Geometry of the arithmetic site, Advances in  Mathematics, {\bf 291}, 2016, 274--329.\\
{[CC2]} A. Connes, C. Consani, Geometry of the scaling site, Selecta Math., (N.S) {\bf 23} (3), 2017, 1803--1850.\\
{[CC3]} A. Connes, C. Consani, Knots, primes and class field theory, to appear  in  Proceedings of the Pisa Conference Regulators V in Contemporary Mathematics, AMS, 2025. \\
{[D1]} C. Deninger, Local $L$-factors of motives and regularized determinants, Inventiones Mathematicae, {\bf 107}, 1992, 135--150.\\
{[D2]} C. Deninger, Lefschetz trace formulas and explicit formulas in analytic number theory, Journal f\"{u}r Reine und Angewandte Mathematik, {\bf 441}, (1993), 1--15.\\
{[D3]} C. Deninger, Motivic $L$-functions and regularized determinants, Proceedings Symposia in Pure Mathematics, AMS, {\bf 55}, 1, 1994, 707--743.\\
{[D4]} C. Deninger, Some analogies between number theory and dynamical 
systems on foliated spaces, Documenta Mathematica, Jahresbericht der Deutschen Mathematiker-Vereinigung, Extra Volume International Congress of Mathematicians I,  
1998, 23--46.\\
{[D5]} C. Deninger, On dynamical systems and their possible significance for
arithmetic geometry, In: Reznikov, A., Schappacher, N. (eds), Regulators in analysis, geometry and number
theory, Progress in  Mathematics {\bf 171}, Birkh\"{a}user Boston, Boston, MA, 2000, 29--87.\\
{[D6]} C. Deninger, Number theory and dynamical systems on foliated spaces,
Jahresbericht der  Deutschen Mathematiker-Vereinigung, {\bf 103}, no. 3, 2001, 79--100.\\
{[D7]} C. Deninger, Dynamical systems for arithmetic schemes, to appear in Indagationes Mathematicae, 2024.\\
{[KS]} R. A. Kucharczyk, P. Scholze, Topological realization of absolute Galois group,  In: Cohomology of arithmetic groups, On the Occasion of Joachim Schwermer's 66th Birthday, Bonn, Germany, June 2016, Springer Proceedings of  Mathematics and  Statistics {\bf 245}, 2018, 201--288.\\
{[M]} M. Morishita, Knots and Primes -- An introduction to Arithmetic Topology, Second Edition,  Universitext. Springer, 2024. \\
{[S]} C.  Soul\'{e}, Sur les z\'{e}ros des fonctions $L$ automorphes, Comptes Rendus de l'Acad\'{e}mie des Sci.ences, Paris Ser. I Math.,  {\bf 328}, 1999, 955--958.
\end{flushleft}
\vspace{0.4cm}
{\small 
M. Morishita:\\
Graduate School of Mathematics, Kyushu University, \\
744, Motooka, Nishi-ku, Fukuoka  819-0395, Japan.\\
e-mail: morishita.masanori.259@m.kyushu-u.ac.jp \\
}

\end{document}